\documentclass[12pt]{amsart}
\usepackage{geometry}
\usepackage{graphicx, amsfonts, amssymb}
\usepackage{mathrsfs, amsmath, amsthm}
\usepackage{verbatim}
\usepackage{relsize}
\usepackage{epsfig}

\def\a{\alpha}

\newcommand{\hol}{\rm{Hol}}
\DeclareMathOperator{\og}{O}

\newcommand{\defeq}{\stackrel{\text{def}}{=}}
\def\D{{\mathbb D}}

\def\C{{\mathbb C}}

\def\Dpa{{\mathcal D^p_{\alpha}}}

\def\({\left(}       \def\){\right)}

\DeclareMathOperator{\op}{o} 

\newcommand{\ig}{\stackrel{\text{def}}{=}}

\newcommand{\n}[1]{\|#1\|}
\newtheorem{theorem}{Theorem}
\newtheorem{lemma}{Lemma}

\newtheorem{proposition}{Proposition}
\theoremstyle{definition}

\newtheorem{remark}{Remark}
\numberwithin{equation}{section}
\theoremstyle{theorem}
\newtheorem{other}{\bf Theorem}              
\newtheorem{otherp}{\bf Proposition}  

\newenvironment{pf}{\noindent{\emph{Proof.}}}{$\Box$ }
\newenvironment{Pf}{\noindent{\emph{Proof of}}}{$\Box$ }
 \DeclareMathOperator{\ogr}{O}

\begin{document}
\title[Rhaly operators acting on Hardy, Bergman, and Dirichlet spaces]
{Rhaly operators acting on Hardy, Bergman, and Dirichlet spaces}


\author[P.~Galanopoulos]{Petros Galanopoulos}
 \address{Department of Mathematics,
Aristotle University of Thessaloniki, 54124, Thessaloniki, Greece}
\email{petrosgala@math.auth.gr}
\author[D.~Girela]{Daniel Girela$^\ast$}
 \address{An\'alisis Matem\'atico,
Universidad de M\'alaga, Campus de Teatinos, 29071 M\'alaga, Spain}
 \email{girela@uma.es}

\subjclass[2020]{Primary 30H10; 47B91; 42B30}
\keywords{Hardy spaces, Bergman spaces, Dirichlet spaces, The
Ces\`{a}ro operator, Rhaly operators, Mean Lipschitz spaces}
\begin{abstract}
In this article we address the question of characterizing the
sequences of complex numbers $(\eta )=\{ \eta_n\}_{n=0}^\infty $
whose associated Rhaly operator $\mathcal R_{(\eta )}$ is bounded or
compact on the Hardy spaces $H^p$ ($1\le p<\infty $), on the Bergman
spaces $A^p_\alpha $, and on the Dirichlet spaces $\mathcal
D^p_\alpha $ ($1\le p<\infty $, $\alpha >-1$). We give a number of
conditions which are either necessary or sufficient for the
boundedness (compactness) of $\mathcal R_{(\eta )}$ on these spaces.
These conditions have to do with the membership in certain mean
Lipschitz spaces of analytic functions of the function $F_{(\eta )}$
defined by $F_{(\eta )}(z)=\sum_{n=0}^\infty \eta _nz^n$ ($z\in
\mathbb D$).
\par
We prove that if $2\le p<\infty $ and $\eta_n=\og \left
(\frac{1}{n}\right )$, then $\mathcal R_{(\eta )}$ is bounded on
$H^p$. However, there exists a sequence $(\eta )$ with $\eta_n=\og
\left (\frac{1}{n}\right )$ such that the operator $\mathcal
R_{(\eta )}$ is not bounded on $H^p$ for $1\le p<2$.
\par We deal also with the derivative-Hardy spaces. For $p>0$ the derivative-Hardy space $S^p$
consists of those functions $f$, analytic in the unit disc $\mathbb
D$, such that $f^\prime \in H^p$. We prove that if $1\le p<\infty $
and $1<q<\infty $ then $\mathcal R_{(\eta )}$ is a bounded operator
from $S^p$ into $S^q$ if and only if it is compact and this happens
if and only if $F_{(\eta )}\in S^q$.
\end{abstract}
\thanks{This research has been supported in part by a grant from \lq\lq El Ministerio de
Ciencia e Innovaci\'{o}n\rq\rq , Spain Project PID2022-133619NB-I00
and by grant from la Junta de Andaluc\'{\i}a FQM-210-G-FEDER)}
\thanks{$^\ast$ Corresponding author: Daniel Girela, girela@uma.es} \maketitle


\section{Introduction and main results}\label{intro}

If $(\eta )=\{ \eta_n\}_{n=0}^\infty $ is a sequence of complex
numbers, the {\it Rhaly matrix} $\mathcal R_{(\eta )}=\mathcal R_{\{
\eta_n\}}$ is defined by
$$\mathcal R_{(\eta )}=\mathcal R_{\{ \eta_n\}}\,=\,\left(\begin{array}{ccccccc}
            \eta_0 & 0  & 0  & 0 & . & . \\
  \eta_1 & \eta_1  & 0  & 0 & . & . \\
  \eta_2 & \eta_2  & \eta_2  & 0 & . & . \\
  \eta_3  &  \eta_3 &  \eta_3 &  \eta_3 & . & . \\
   .  & . & . & . & . & . \\.  & . & . & . & . & . \\
\end{array}\right).$$
This matrix induces the Rhaly operator $\mathcal R_{(\eta
)}=\mathcal R_{\{ \eta_n\}}$ by matrix multiplication, that is, if
$(a)=\{ a_n\}_{n=0}^\infty $ is a sequence of complex numbers
$$\mathcal R_{(\eta )}(a)\,=\,\mathcal R_{ \{\eta_n\} }(a)\,=\,\left \{
\eta_n\sum_{k=0}^na_k\right \}_{n=0}^\infty .$$ The case when
$\eta_n=\frac{1}{n+1}$ for all $n$ corresponds to the classical
Ces\`{a}ro operator $\mathcal C$ which is known to be bounded from
$\ell^p$ to $\ell^p$ for $1<p\le \infty $. This was proved by Hardy
\cite{H-20} and Landau \cite{L} (see also \cite[Theorem\,\@326,
p.\,\@239~]{HLP}). The sequences $\{ \eta_n\}_{n=0}^\infty $ for
which the operator $\mathcal R_{\{ \eta _n\} }$ is either bounded or
compact on the space $\ell^p$ for any $p\in (1, \infty )$ have been
characterized in \cite{GGPr}. A characterization of the sequences
$\{ \eta_n\}$ for which $\mathcal R_{\{ \eta _n\} }$ belongs to the
Schatten class $\mathcal S^q(\ell^2)$ ($1<q<\infty $) has been
obtained in \cite{BDS}.

\par\medskip
Let $\D=\{z\in\C: |z|<1\}$ denote the open unit disc in the complex
plane $\C$. If $\mu $ is a complex Borel measure in $\mathbb D$ and,
for $n=0, 1, 2, \dots $, $\mu_n$ denotes the moment of order $n$ of
$\mu $, that is,
$$\mu_n\,=\,\int_{\mathbb D}w^n\,d\mu (w),\quad n=0, 1, 2, \dots ,$$
the operator $\mathcal R_{\{\mu _n\} }$ will be denoted by $\mathcal
C_{\mu }$. When $\mu $ is the Lebesgue measure on the radius $[0,
1)$ the operator $\mathcal C_{\mu }$ is the Ces\`{a}ro operator
$\mathcal C$.
\par\medskip
Let $\hol (\D)$ be the space of all analytic functions in $\D$
endowed with the topology of uniform convergence in compact subsets.
\par\medskip If $\,0<r<1\,$ and $\,f\in \hol (\D)$, we set
$$
M_p(r,f)=\left(\frac{1}{2\pi }\int_0^{2\pi }
|f(re^{it})|^p\,dt\right)^{1/p}, \,\,\, 0<p<\infty ,
$$
$$
M_\infty(r,f)=\sup_{\vert z\vert =r}|f(z)|.
$$
For $\,0<p\le \infty $,\, the Hardy space $H^p$ consists of those
$f\in \hol(\mathbb D)$ such that $$\Vert f\Vert _{H^p}\defeq
\sup_{0<r<1}M_p(r,f)<\infty .$$ We refer to \cite{Du:Hp} for the
notation and results regarding Hardy spaces.
\par\medskip
Let $dA$ denote the area measure on $\mathbb D$, normalized so that
the area of $\mathbb D$ is $1$. Thus $dA(z)\,=\,\frac{1}{\pi
}\,dx\,dy\,=\,\frac{1}{\pi }\,r\,dr\,d\theta $. For $0<p<\infty $
and $\alpha >-1$ the weighted Bergman space $A^p_\alpha $ consists
of those $f\in \hol (\mathbb D)$ such that
$$\Vert f\Vert _{A^p_\alpha }\,\defeq\, \left (\int_\mathbb D\vert f
(z)\vert ^p\,dA_\alpha (z)\right )^{1/p}\,<\,\infty ,$$ where
$dA_\alpha (z)=(\alpha +1)(1-\vert z\vert ^2)^{\alpha }dA(z)$. We
refer to \cite{DS,HKZ,Zhu} for the notation and results about
Bergman spaces.
\par
If $0<p<\infty $ and $\a>-1$ the space of Dirichlet type $\Dpa$
consists of all primitives of functions in $A^p_\alpha$.  Hence, if
$f$ is analytic in $\mathbb D$, then $f\in \Dpa$ if and only if
$$
\n{f}_{\Dpa}^p\ig |f(0)|^p+\n{f'}_{A^p_\alpha }^p <\infty .
$$
\par
\par\medskip Identifying any given function $f\in \hol (\mathbb D)$
with the sequence $\{ a_k\}_{k=0}^\infty $ of its Taylor
coefficients, the Ces\`{a}ro operator $\mathcal C$ and the operators
$\mathcal C_\mu $ associated to complex Borel measures on $\mathbb
D$ become linear operators from $\hol (\mathbb D)$ into itself as
follows:
\par If $f\in \hol (\mathbb D)$, $f(z)=\sum_{k=0}^\infty a_kz^k$
\,($z\in \mathbb D$), then
$$\mathcal C(f)(z)=\sum_{n=0}^\infty \left
(\frac{1}{n+1}\sum_{k=0}^na_k\right )z^n,\quad \mathcal C_\mu
(f)(z)=\sum_{n=0}^\infty \left (\mu _n\sum_{k=0}^na_k\right
)z^n,\quad z\in \mathbb D.$$ The Ces\`{a}ro operator is known to be
bounded on the Hardy spaces $H^p$ ($0<p<\infty $) and on the Bergman
spaces $A^p_\alpha $ ($0<p<\infty $, $\alpha >-1$). This has been
established by a number of authors using different methods
\cite{Andersen, Mi, Nowak, Sis-1987, Sis-1990, Sis-1996, Stem}.
\par The operators $\mathcal C_\mu $ with $\mu $ being a finite
positive Borel measure supported on the radius $[0,1)$ were
introduced in \cite{GGM-2022} where it was proved that such a
$\mathcal C_\mu $ is bounded on the Hardy space $H^p$ ($1\le p<
\infty $), or on the Bergman space $A^p_\alpha $ ($1<p<\infty $,
$\alpha >-1$) if and only if $\mu $ is a Carleson measure,  that is,
if there exists a positive constant $C$ such that $\mu ([r, 1))\,\le
\, C(1-r)$ for $0<r<1$ or, equivalently, if $\mu_n=\og \left
(\frac{1}{n+1}\right )$.  Theorem~9 of \cite{GG-2025} shows that
this is also true for the spaces $A^1_\alpha $ ($\alpha >-1$). Since
its introduction, the boundedness of these operators $\mathcal C_\mu
$ on other spaces of analytic functions in $\mathbb D$ has been
studied by distinct authors (see, e.\,\@g., \cite{BaoSWu, BBJ, GGMM,
JT, T}).
\par Blasco \cite{Blasco-2024} considered the operators $\mathcal
C_\mu $, with $\mu $ being a complex Borel measure on $[0, 1)$,
acting on Hardy spaces. Later on, a good amount of work has been
carried out considering the operators $\mathcal C_\mu $, where $\mu
$ is a complex Borel measure  on $\mathbb D$, acting on distinct
subspaces of $\hol (\mathbb D)$; see, for instance,
\cite{Blasco-2024-1, GGM-2023, LX}.
\par\medskip Let us remark that when $\mu $ is a complex Borel
measure on $\mathbb D$, the operator $\mathcal C_\mu $ has an
integral representation: If $f\in \hol (\mathbb D)$,
$f(z)=\sum_{k=0}^\infty a_kz^k$ ($z\in \mathbb D$), then
\begin{equation}\label{intrep}\mathcal C_\mu (f)(z)\,=\,\sum_{n=0}^\infty
\mu_n\left( \sum_{k=0}^na_k\right )z^n\,=\,\int_{\mathbb
D}\frac{f(wz)}{1-wz}\,d\mu (w),\quad z\in \mathbb D.\end{equation}
This plays a crucial role to obtain a good number of the results we
have just mentioned.
\par\medskip
If $(\eta )=\{ \eta _n\}_{n=0}^\infty $ is a sequence of complex
numbers then the Rhaly operator $\mathcal R_{\{ \eta_n\}}$ can be
viewed  formally as an operator acting on spaces  of analytic
functions in the disc: If $f\in \hol (\mathbb D)$, $f(z)=\sum
_{n=0}^\infty a_nz^n $ ($z\in \mathbb D$), formally, we set
$$\mathcal R_{(\eta )}(f)(z)\,=\,\mathcal R_{\{ \eta _n\}}(f)(z)=\sum_{n=0}^\infty \eta _n\left
(\sum_{k=0}^na_k\right )z^n,\quad z\in\mathbb D,$$ whenever the
right hand side makes sense and defines an analytic function in
$\mathbb D$. \par The sequence $(\eta )=\{ \eta _n\}$ may not be
bounded. Hence, in contrary to what happens with the operators
$\mathcal C_\mu $, we cannot assure that the operator $\mathcal
R_{\{ \eta_n\}}$ is a well defined operator from $\hol (\mathbb D)$
into itself. Proposition~1 of \cite{Bl-Ga-Gi} asserts that if $X$ is
a subspace of $\hol (\mathbb D)$ which contains the constants, then
$\mathcal R_{\{ \eta_n\}}$ is a well defined linear operator from
$X$ into $\hol (\mathbb D)$ if and only if the power series
$\sum_{n=0}^\infty \eta _nz^n$ defines an analytic function in
$\mathbb D$.
\par From now on, if the power series $\sum_{n=0}^\infty \eta _nz^n$ defines an
analytic function in $\mathbb D$, we shall set
$$F_{(\eta )}(z)\,=\,F_{\{\eta _n\}}(z)\,=\,\sum_{n=0}^\infty \eta
_nz^n,\quad z\in \mathbb D.$$
\par\medskip In view of the identification of $H^2$ and $\ell^2$,
the sequences $(\eta )$ for which $\mathcal R_{(\eta )}$ is bounded
or compact on $H^2$ are characterized in \cite{GGPr}. The articles
\cite{BaoGSWa, Bl-Ga-Gi, GG-2025} contain characterizations of the
sequences $(\eta )$ for which the operator $\mathcal R_{(\eta )}$ is
bounded (compact) from  $\mathcal D^2_\alpha $ into $\mathcal
D^2_\beta $ ($\alpha , \beta \in \mathbb R$). Here, for $\alpha\in
\mathbb R$, $\mathcal D^2_\alpha$ is the space of functions
$f\in\hol(\mathbb D)$ such that $|a_0|^2+\sum_{n=1}^\infty
n^{1-\alpha} |a_n|^2<\infty$, where $f(z)=\sum_{n=0}^\infty a_nz^n$.
\par
When we deal with a general Rhaly operator $\mathcal R_{(\eta )}$,
we do not have an integral representation similar to (\ref{intrep}).
This is one reason why a characterization of those sequences $(\eta
)$ for which $\mathcal R_{(\eta )}$ is bounded (compact) on $H^p$ or
on $A^p_\alpha $ with $p\neq 2$ has not been obtained yet. In this
paper we address the question of characterizing such sequences for
$1<p<\infty $. For simplicity, if $X$ and $Y$ are two Banach spaces,
$\mathcal B(X,Y)$ will denote the space of continuous linear
operators from $X$ into $Y$ and $\mathcal K(X,Y)$ will stand for the
space of all compact linear operators from $X$ into $Y$. Also
$\mathcal B(X)$ will stand for $\mathcal B(X,X)$ and $\mathcal K(X)$
will stand for $\mathcal K(X,X)$. Our main results are the
following.

\begin{theorem}\label{boundedness-Hp} Let $(\eta )=\{ \eta _n\}_{n=0}^\infty $ be a
sequence of complex numbers. \par (a) If $1<p<\infty $ and $\mathcal
R_{(\eta )}\in \mathcal B(H^p)$ then
\begin{equation}\label{condLomega}M_p(r, F_{(\eta )}^\prime )\,=\,\og \left
(\frac{\log\frac{1}{1-r}}{(1-r)^{1-\frac{1}{p}}}\right ),\quad
\text{ as $r\to 1$.}\end{equation}
\par (b) If $1<p\le 2$ and $F_{(\eta)} \in \Lambda^p_{1/p}$ then
$\mathcal R_{(\eta )}\in \mathcal B(H^p)$.
\par (c) If $2<p<\infty $ and $F_{(\eta )} \in \Lambda ^q_{1/q}$ for some
$q\in (2, p)$, then $\mathcal R_{(\eta )}\in \mathcal B(H^p)$.
\par (d) $\mathcal R_{(\eta )}\in \mathcal B(H^2)$ if and only if
$F_{(\eta )}\in \Lambda^2_{1/2}$.
\end{theorem}
\par\medskip
\begin{theorem}\label{compactnessHp} Let $(\eta )=\{ \eta _n\}_{n=0}^\infty $ be a
sequence of complex numbers. \par (a) If $1<p<\infty $ and $\mathcal
R_{(\eta )}\in \mathcal K(H^p)$ then
\begin{equation}\label{condlomega}M_p(r, F_{(\eta )}^\prime )\,=\,\op \left
(\frac{\log\frac{1}{1-r}}{(1-r)^{1-\frac{1}{p}}}\right ),\quad
\text{ as $r\to 1$.}\end{equation}
\par (b) If $1<p\le 2$ and $F_{(\eta)} \in \lambda^p_{1/p}$ then
$\mathcal R_{(\eta )}\in \mathcal K(H^p)$.
\par (c) If $2<p<\infty $ and $F_{(\eta )} \in \lambda ^q_{1/q}$ for some
$q\in (2, p)$, then $\mathcal R_{(\eta )}\in \mathcal
K(H^p)$.
\par (d) $\mathcal R_{(\eta )}\in \mathcal K(H^2)$ if and only if
$F_{(\eta )}\in \lambda^2_{1/2}$.\end{theorem}
\par\medskip
\begin{theorem}\label{Apa} Let $(\eta )\,=\,\{
\eta_n\}_{n=0}^\infty $ be a sequence of complex numbers.
\par (i) If $\alpha >-1$, $1\le p<\infty $, and $F_{(\eta )}\in
\Lambda ^p_{1/p}$ then $\mathcal R_{(\eta )}\in \mathcal
B(A^p_\alpha )$.
\par (ii) If $\alpha >-1$, $1<p<\infty $, and $F_{(\eta )}\in
\lambda ^p_{1/p}$ then $\mathcal R_{(\eta )}\in \mathcal
K(A^p_\alpha )$.
\par (iii) If $1<p<\infty $, $-1<\alpha <2p-2$, and $\mathcal R_{(\eta )}\in
\mathcal B(A^p_\alpha )$ then (\ref{condLomega}) holds.
\par (iv) If $1<p<\infty $, $-1<\alpha <2p-2$, and $\mathcal R_{(\eta )}\in
\mathcal K(A^p_\alpha )$ then (\ref{condlomega}) holds.
\par (v) If $\alpha >-1$, then $\mathcal R_{(\eta )}\in \mathcal B(A^2_\alpha )$ if and only if
$F_{(\eta )}\in \Lambda^2_{1/2}$.
\par (vi) If $\alpha >-1$, then $\mathcal R_{(\eta )}\in \mathcal K(A^2_\alpha )$ if and only if
$F_{(\eta )}\in \lambda^2_{1/2}$.
\end{theorem}
\par\medskip
\begin{theorem}\label{Dpa} Let $(\eta )\,=\,\{
\eta_n\}_{n=0}^\infty $ be a sequence of complex numbers.
\par (i) If $1\le p<\infty $, $\alpha >p-2$,  and $F_{(\eta )}\in
\Lambda ^p_{1/p}$ then $\mathcal R_{(\eta )}\in \mathcal B(\mathcal
D^p_\alpha )$.
\par (ii) If $1\le p<\infty $, $\alpha >p-2$, and $F_{(\eta )}\in
\lambda ^p_{1/p}$ then $\mathcal R_{(\eta )}\in \mathcal K(\mathcal
D^p_\alpha )$.
\par (iii) If $1<p<\infty $, $-1<\alpha <2p-2$, and $\mathcal R_{(\eta )}\in
\mathcal B(\mathcal D^p_\alpha )$ then (\ref{condLomega}) holds.
\par (iv) If $1<p<\infty $, $-1<\alpha <2p-2$, and $\mathcal R_{(\eta )}\in
\mathcal K(\mathcal D^p_\alpha )$ then (\ref{condlomega}) holds.
\par (v) If $\alpha >0$, then $\mathcal R_{(\eta )}\in \mathcal B(\mathcal D^2_\alpha )$ if and only if
$F_{(\eta )}\in \Lambda^2_{1/2}$.
\par (vi) If $\alpha >0$, then $\mathcal R_{(\eta )}\in \mathcal K(\mathcal D^2_\alpha )$ if and only if
$F_{(\eta )}\in \lambda^2_{1/2}$.
\end{theorem}
\medskip
Since the Ces\`{a}ro operator is bounded on $H^p$ for\, $0<p<\infty
$, it is natural to ask whether or not the condition $\eta
_n\,=\,\og \left (\frac{1}{n}\right )$ implies that $\mathcal
R_{(\eta )}\in \mathcal B(H^p)$ for some (all) $p$. We shall prove
the following result.
\begin{theorem}\label{og1n}
\par (i) Suppose that $2\le p<\infty $ and let $(\eta )=\{ \eta _n\}_{n=0}^\infty $ be a sequence of
complex numbers. If\, $\eta_n\,=\,\og \left (\frac{1}{n}\right )$,
as $n\to \infty $, then $\mathcal R_{(\eta )}\in \mathcal B(H^p)$.
\par (ii) There exists a sequence of complex numbers $(\eta
)=\{ \eta _n\}_{n=0}^\infty $,  with $\eta_n\,=\,\og \left
(\frac{1}{n}\right )$, as $n\to \infty $, such that $\mathcal
R_{(\eta )}\not\in \mathcal B(H^p)$ for $1\le p<2$.

\par (iii) If\, $1<p<\infty $\, and $(\eta )=\{ \eta _n\}_{n=0}^\infty
$ is a decreasing sequence of non negative real numbers then
$\mathcal R_{(\eta )}\in \mathcal B(H^p)$ if and only if
$\eta_n\,=\,\og \left (\frac{1}{n}\right )$, as $n\to \infty $.
\end{theorem}
\par\medskip Finally, we are going to deal with the derivative-Hardy
spaces $S^p$. For $p>0$, the derivative-Hardy space $S^p$ consists
of those $f\in \hol (\mathbb D)$ such that $f^\prime \in H^p$. For
$f\in S^p$, we set
$$\Vert f\Vert _{S^p}^p\,=\,\vert f(0)\vert ^p\,+\,\Vert f^\prime
\Vert _{H^p}^p.$$ In some instances the space $S^p$ plays the role
of \lq\lq the limit of $\mathcal D^p_\alpha $, as $\alpha \to
-1^+$\rq\rq . Corollary\,\@1 of \cite{Bl-Ga-Gi} asserts that
$$\mathcal R_{(\eta)}\in \mathcal B(S^2)
\,\,\, \Leftrightarrow\,\,\, \mathcal R_{(\eta)}\in \mathcal
K(S^2)\,\,\, \Leftrightarrow\,\,\, F_{(\eta )}\in S^2.$$ We extend
this in the following theorem.
\begin{theorem}\label{RSp} Suppose that $1\le p<\infty $,
$1<q<\infty $, and let $(\eta )=\{ \eta _n\}_{n=0}^\infty $ be a
sequence of complex numbers. Then the following conditions are
equivalent.
\par (i) $\mathcal R_{(\eta)}\in \mathcal B(S^{p}, S^{q})$.
\par (ii) $\mathcal R_{(\eta)}\in \mathcal K(S^{p}, S^{q})$.
\par (iii) $F_{(\eta )}\in S^{q}$.
\end{theorem}
\par\medskip The mean Lipschitz spaces $\Lambda^p_{1/p}$, the \lq\lq little oh\rq\rq mean Lipschitz spaces
$\lambda^p_{1/p}$, and, more generally the mean Lipschitz spaces
$\Lambda (p,\omega )$ will be properly defined in
Section\,\@\ref{PDR}, which will be devoted to introduce a number of
spaces and results which will be needed in our work. The proofs of
the results for $p=2$ will be presented in Section\,\@\ref{pequal2}.
Section\,\@\ref{necessary} will include the proofs of the necessity
statements in Theorems \ref{boundedness-Hp},
 \ref{compactnessHp}, \ref{Apa}, and \ref{Dpa} for $p\neq 2$. The sufficiency
 statements in these theorems will be proved in Section\,\@\ref{sufficient}.
 Theorem\,\@\ref{og1n} will be proved in Section\,\@\ref{further}
where we will present also a number of results regarding the
operators $\mathcal R_{(\eta )}$ acting on $H^1$. Finally, the proof
of Theorem\,\@\ref{RSp} will be presented in
Section\,\@\ref{RhalySp}.
\par
 Throughout the paper
 we shall be using the convention that
$C=C(p, \alpha ,q,\beta , \dots )$ will denote a positive constant
which depends only upon the displayed parameters $p, \alpha , q,
\beta \dots $ (which sometimes will be omitted) but not  necessarily
the same at different occurrences. Furthermore, for two real-valued
functions $K_1, K_2$ we write $K_1\lesssim K_2$, or $K_1\gtrsim
K_2$, if there exists a positive constant $C$ independent of the
arguments such that $K_1\leq C K_2$, respectively $K_1\ge C K_2$. If
we have $K_1\lesssim K_2$ and $K_1\gtrsim K_2$ simultaneously, then
we say that $K_1$ and $K_2$ are equivalent and we write $K_1\asymp
K_2$.
\par\medskip

\section{Preliminary definitions and results}\label{PDR} \par In
this section we are going to present a number of spaces and results
which will be needed to prove the theorems stated in
section\,\@\ref{intro}. Let start fixing some notation. If $f\in
\hol (\mathbb D)$, $f(z)=\sum_{k=0}^\infty a_kz^k$ ($z\in \mathbb
D)$, we set
\begin{equation}\label{deltanm}S_{n,m}f(z)\,=\,\sum_{k=n}^ma_kz^k,\quad
0\le n<m<\infty .\end{equation}
\begin{equation}\label{deltaN}\Delta _Nf(z)\,=\,S
_{N,2N-1}f(z)\,=\,\sum_{k=N}^{2N-1}a_kz^k, \quad 1\le N<\infty
.\end{equation}
\begin{equation}\label{SNRN}S _Nf(z)\,=\,\sum_{k=0}^{N}a_kz^k, \quad 0\le N<\infty .\end{equation}
\subsection{Mean Lipschitz spaces}\label{Mean Lip}
If $f$ is a function which is analytic in $\mathbb D $ and has a
non-tangential limit $f(e\sp {i\theta })$ at almost every $\xi \in
\partial \mathbb D$, we define
$$
\aligned \omega _ p(\delta , f)= \sup_ {0<\vert t\vert \leq \delta
}\left (\frac{1}{2\pi }\int _{-\pi }\sp \pi \left \vert f(e\sp
{i(\theta +t)})-f(e\sp {i\theta })\right \vert \sp p\,
d\theta \right )\sp {1/p}, &\quad \delta >0,\quad\hbox{if $1\leq p<\infty $},\\
\omega _ \infty (\delta , f)= \sup_ {0<\vert t\vert \leq \delta
}\left (\operatornamewithlimits{ess.sup}_ {\theta\in [-\pi, \pi ]}
\vert f(e\sp {i(\theta +t)})-f(e\sp {i\theta })\vert\right ), &\quad
\delta >0.\endaligned $$ Then $\omega _ p(., f)$ is the integral
modulus of continuity of order $p$ of the boundary values $f(e\sp
{i\theta })$ of $f$.\par Throughout the paper $\omega :[0,
1]\rightarrow [0, \infty )$ will be a continuous and increasing
function with $\omega (0)=0$. Then, for $1\leq p\leq \infty $, the
mean Lipschitz space $\Lambda (p, \omega )$ consists of those
functions $f\in H\sp p$ which satisfy
$$\omega_ p(\delta , f)=\og (\omega (\delta )),\quad\hbox{as $\delta \to 0$}.$$
If $0<\alpha \le 1$ and $\omega (\delta )=\delta \sp\alpha $, we
shall write $\Lambda \sp p_ \alpha $ instead of $\Lambda (p, \omega
)$, that is, we set
$$\Lambda \sp p_ \alpha =\Lambda (p, \delta\sp\alpha ),\quad 0<\alpha \le
1,\quad 1\le p\le \infty .$$ The corresponding \lq\lq little
oh\rq\rq \, spaces are denoted by $\lambda ^p_\alpha $ and $\lambda
(p, \omega)$. \par If $0<\alpha \le 1$, $1\le p\le \infty $, and
$f\in \hol (\mathbb D)$, we set
$$\beta _{p,\alpha }(r, f)\,=\,(1-r)^{1-\alpha }M_p(r, f^\prime
),\quad 0<r<1.$$ \par
 Classical results of Hardy and Littlewood (see \cite{BSS} and
\cite[Chapter\,\@5]{Du:Hp}) show that whenever $0<\alpha \le 1$ and
$1\le p\le \infty $ we have that $\Lambda ^p_{\alpha }\subset H^p$
and that
\begin{equation}\Lambda ^p_{\alpha }\,=\,\left \{ \hbox{$f$ analytic
in $\D $: $M_ p(r, f\sp\prime )= \og \left ((1-r)\sp {\alpha
-1}\right )$, \quad as $r\to 1$} \right\}.\end{equation} This result
was extended in \cite{BS} where Blasco and de Souza proved that if
$\omega $ satisfies the so called Dini and $b_1$ conditions then
$$\Lambda (p, \omega )=\left \{ \hbox{$f$ analytic in $\Delta $ :
$M_ p(r, f\sp\prime )=\og \left ( \frac{\omega (1-r)}{1-r}\right ),$
as $r\to 1$}\right \}.$$ In particular this is true for $\omega $
defined by
\begin{equation}\label{omegaalpha}\omega (\delta )\,=\,\delta ^{\alpha
}\log\frac{A}{\delta },\quad 0<\delta <1,\end{equation} with
$0<\alpha <1$ and $A\ge e^{1/\alpha }$. Thus, if $\omega $ is as in
(\ref{omegaalpha}) with $0<\alpha <1$ and $A\ge e^{1/\alpha }$, then
\begin{equation}\label{lampomega}\Lambda (p, \omega )\,=\,\left \{ f\in H^p: M_p(r,
f^\prime )\,=\,\og \left (\frac{\log\frac{1}{1-r}}{(1-r)^{1-\alpha
}}\right ),\,\,\text{as $r\to 1$}\right \} .\end{equation}
\par If $1<p<\infty $ and $(1/p)<\alpha \le 1$, then
$\Lambda^p_\alpha \subset \Lambda_{\alpha -{(1/p)}}$ and, hence,
each function in $\Lambda ^p_\alpha $ has a continuous extension to
the closed unit disc \cite[p.\,\@88]{BSS}. This is not true for the
space $\Lambda ^p_{1/p}$. This follows easily noticing that the
function $f(z)=\log (1-z)$ is an unbounded function which belongs to
$\Lambda_ {1/p}\sp p$ for all $p\in (1, \infty )$. Bourdon, Shapiro,
and Sledd \cite{BSS} proved that $\Lambda^p_{1/p}\subset BMOA$ for
all $p\in (1, \infty )$. This inclusion was shown to be sharp in a
very strong sense in \cite{BGM}. Let us remark also that the
$\Lambda^p_{1/p}$-spaces form a nested scale of spaces
\begin{equation}\label{nested}
\Lambda ^q_{1/q}\,\subset \,\Lambda ^p_{1/p}\quad 1\le q<p<\infty .
\end{equation}
\par\medskip
For $1\le p\le \infty $ and $0<\alpha \le 1$, the space $\Lambda
^p_\alpha $ is a Banach space with the norm $\Vert \cdot \Vert _{p,
\alpha }$ defined by
$$\Vert f\Vert _{p, \alpha }\,=\,\vert f(0)\vert \,+\,\beta_{p,
\alpha }(f),$$ where
$$\beta_{p,
\alpha }(f)\,=\,\sup_{0<r<1}(1-r)^{1-\alpha }M_p(r, f^\prime
)\,=\,\sup_{0<r<1}\beta_{p,\alpha }(r, f).$$ Let us recall several
distinct characterizations of $\Lambda^p_\alpha$ spaces, (see
\cite{BSS}, \cite{Du:Hp}, \cite{GGPS}, and \cite{MP}).

\begin{other}\label{th:mlip}
Suppose that $1<p<\infty$, $0<\alpha<1$ and $g\in \hol(\D)$. The
following conditions are equivalent:
\begin{enumerate}
\item[(i)] $g\in \Lambda^p_\alpha .$
\item[(ii)] $M_p(r,g')=\ogr\left(\frac{1}{(1-r)^{1-\alpha}}\right)$, as $r\to 1^-$.
\item[(iii)] $\|\Delta_Ng\|_{H^p}=\ogr\left(N^{-\alpha}\right)$, as $N\to\infty$.
\item[(iv)] $\|\Delta_Ng'\|_{H^p}=\ogr\left(N^{(1-\alpha)}\right)$, as $N\to\infty$.
\end{enumerate}
In fact,
$$\beta _{p, \alpha }(f)\,\asymp \,\sup_{N\ge 1}N^\alpha \|\Delta_Ng\|_{H^p}\,\asymp \,\sup_{N\ge
1}N^{\alpha -1}\|\Delta_Ng'\|_{H^p}.$$
\end{other}

\begin{remark}\label{re:mlip}
The corresponding results for the little-oh space $\lambda^p_\alpha
$ remain true, and they can be proved following the proofs in the
references for Theorem \ref{th:mlip}.
\end{remark}

The following result will be needed in our work. Since we have not
found any reference to it in the literature, we will include a
proof.
\begin{proposition}\label{partialsumlamba} Suppose that $1<p<\infty
$, $0<\alpha \le 1$, and $f\in \lambda ^p_\alpha $. Then:
\par (i) Defining, for $0<r<1$, $f_r$ by
$f_r(z)=f(rz)$\, ($\vert z\vert \le 1$), we have that
$$\beta_{p, \alpha }(f-f_r)\,\rightarrow 0,\quad \text{as $r\to
1$}.$$
\par (ii)
$S_Nf\,\rightarrow f,\,\,\text{as $N\to \infty $, in the norm of
$\Lambda^p_\alpha $}.$
\end{proposition}
\begin{pf} Let us start proving (i). Take  a non-constant $f\in \lambda^p_\alpha $. Take $\varepsilon >0$
with $0<\varepsilon <3\beta_{p, \alpha
}(f).$ Since $f\in \lambda ^p_\alpha $, there exists $s_0\in (0, 1)$
such that \begin{equation}\label{s0} (1-s)^{1-\alpha }M_p(s,
f')\,<\frac{\varepsilon }{3}, \quad s_0\le s<1.\end{equation} For
$r, s\in (0,1)$, we have
\begin{align}\label{first}&(1-s)^{1-\alpha }M_p(s, f'-f'_r)\,=\,(1-s)^{1-\alpha
}\left [\frac{1}{2\pi}\int_0^{2\pi }\left \vert f'(se^{i\theta
})-rf'(rse^{i\theta })\right \vert ^p\,d\theta \right ]^{1/p}\\
=\,&(1-s)^{1-\alpha }\left [\frac{1}{2\pi}\int_0^{2\pi }\left \vert
f'(se^{i\theta })-rf'(se^{i\theta })+rf'(se^{i\theta
})-rf'(rse^{i\theta })\right \vert ^p\,d\theta \right
]^{1/p}\nonumber \\ \le \,&(1-r)(1-s)^{1-\alpha }M_p(s,
f')\,+\,r(1-s)^{1-\alpha }\left [\frac{1}{2\pi }\int_0^{2\pi }\left
\vert f'(se^{i\theta })-f'(rse^{i\theta })\right \vert ^p\,d\theta
\right ]^{1/p}.\nonumber
\end{align}
Using this, (\ref{s0}), and the fact that $M_p(rs, f')\le M_p(s,
f')$ ($0<s, r<1$), we obtain
\begin{align}\label{s0s1}(1-s)^{1-\alpha }M_p(s, f'-f'_r)\,\le \,&(1-s)^{1-\alpha
}M_p(s, f')+\,(1-s)^{1-\alpha }\left [M_p(s, f')\,+\,M_p(rs,
f')\right ] \nonumber\\\le\,&3(1-s)^{1-\alpha }M_p(s,
f')\\<\,&\varepsilon ,\qquad s_0<s<1.\nonumber\end{align}
\par Since $f'$ is uniformly continuous in the closed disc $\{ z\in
\mathbb D: \vert z\vert \le s_0\} $, there exists $\delta \in (0,
1)$ such that
\begin{equation}\label{unifcont}\vert f'(z)-f'(w)\vert
<\frac{2\varepsilon}{3},\,\,\text{if $\vert z\vert \le s_0, \,\vert
w\vert \le s_0$,\,\,\,and $\vert z-w\vert \le \delta
$}.\end{equation} Take $r_0\in (0, 1)$ with $$r_0>\max \left
(1-\delta , 1-\frac{\varepsilon }{3\beta_{p, \alpha }(f)}\right ).$$
Using (\ref{first}) and (\ref{unifcont}), we see that, if $0<s\le
s_0$ and $r_0<r<1$, then
$$(1-s)^{1-\alpha }M_p(s, f'-f'_r)\le \,(1-r)\beta_{p, \alpha
}(f)\,+\,(1-s)^{1-\alpha }\frac{2\varepsilon }{3}\le
\frac{\varepsilon}{3\beta_{p, \alpha }(f)}\beta_{p, \alpha
}(f)+\,\frac{2\varepsilon }{3}=\,\varepsilon.$$ This and
(\ref{s0s1}) give that $\beta _{p, \alpha }(f-f_r)<\varepsilon $, if
$r_0<r<1$. Hence, (i) is proved.
\par\medskip In order to prove (ii), we shall start proving the
following weaker result. \begin{lemma}\label{seqpol} Suppose that
$1<p<\infty $, $0<\alpha \le 1$, and $g\in \lambda ^p_\alpha $. Then
$g$ is the limit in $\Lambda ^p_\alpha $ of a sequence of
polynomials.\end{lemma}\begin{pf} Take $g\in \lambda ^p_\alpha $.
Set $r_n=1-\frac{1}{n}$ ($n=2, 3, \dots $) and consider the sequence
of functions $\{ g_{r_n}\} _{n=2}^\infty $. For every $n\ge 2$,
$g_{r_n}$ is holomorphic in $\{ \vert z\vert \le \frac{1}{r_n}\} $
and, hence, there exists a polynomial $P_n$ with $P_n(0)=g(0)$ such
that
$$\vert g'_{r_n}(z)\,-\,P'_n(z)\vert \le \frac{1}{n},\quad \vert z\vert
\le 1.$$ Clearly, this implies that
$$\beta_{p, \alpha }(g_{r_n}-P_n)\le \frac{1}{n},\quad n\ge 2.$$
Using this and (i) we see that
$$\beta _{p, \alpha }(g-P_n)\,\le \beta _{p, \alpha
}(g-g_{r_n})\,+\,\beta _{p, \alpha }(g_{r_n}-P_n)\,\to
0,\quad\text{as $n\to \infty $}.$$ Thus, $P_n\,\to g$, as $n\to
\infty $, in the norm of $\Lambda^p_\alpha $.
\end{pf}
\par\medskip Bearing in mind Lemma\,\@\ref{seqpol} and Proposition
\,\@1 of \cite{Zhu-Sn}, (ii) will follow if we prove that there
exists a positive constant $C$ such that
\begin{equation}\label{SnLambda}\beta _{p, \alpha }(S_nf)\le C\beta _{p,
\alpha }(f),\,\text{
 for all
$f\in \Lambda ^p_\alpha $ and all $n\in \mathbb N$.}\end{equation}
\par Since $1<p<\infty $, using the Riesz projection theorem, it
follows that the operators $S_n$ are uniformly bounded on $H^p$,
that is, there exists $A>0$ such that
\begin{equation}\label{SnHp}\Vert S_nf\Vert _{H^p}\le A\Vert f\Vert
_{H^p},\quad (n\in \mathbb N, \,\,f\in H^p).\end{equation} Take now
$f\in \Lambda ^p_\alpha $, $f(z)\,=\,\sum_{k=0}^\infty a_kz^k$
($z\in \mathbb D$). Notice that for $n\in\mathbb N$ we have
$$\left (S_nf\right
)'(z)\,=\,\sum_{k=0}^{n-1}(k+1)a_{k+1}z^k\,=\,S_{n-1}(f')(z),\quad
z\in \mathbb D.$$ Then, using (\ref{SnHp}), we see that
\begin{align*}&\beta_{p,\alpha }(S_nf)\,=\,\sup_{0<r<1}(1-r)^{1-\alpha
}M_p(r, (S_nf)')\\=\,&\sup_{0<r<1}(1-r)^{1-\alpha }M_p(r,
S_{n-1}(f'))\,\le A\sup_{0<r<1}(1-r)^{1-\alpha }M_p(r,
f')\,=\,A\beta_{p, \alpha }(f).\end{align*} Hence, (\ref{SnLambda})
holds with $C=A$. This finishes the proof.
\end{pf}
\par\medskip
\subsection{Convolution of analytic functions}\label{multiplier}
 \par If $f$
and $g$ are two analytic functions in the unit disc,
$$f(z)=\sum_{n=0}^\infty a_nz^n,\quad g(z)=\sum_{n=0}^\infty
b_nz^n,\quad z\in \mathbb D,$$ the convolution $f\star g$ of $f$ and
$g$ is defined by \begin{equation}\label{convdef}f\star
g(z)\,=\,\sum_{n=0}^\infty a_nb_nz^n,\quad z\in \mathbb
D.\end{equation} We have \begin{equation}\label{convdefint}f\star
g(\rho re^{i\theta })\,=\,\frac{1}{2\pi }\int_0^{2\pi }f(\rho
e^{it})g(re^{i(\theta -t)})\,dt,\quad 0<\rho <1,\,\,0<r<1,\,\,\theta
\in \mathbb R.\end{equation} It is well known that the convolution
of a function in $H^1$ and another one in $H^p$ ($p\ge 1$) lies in
$H^p$ and that
\begin{equation*}\label{conv}\Vert f\star g\Vert_{H^p}\le \Vert
f\Vert _{H^1}\Vert g\Vert _{H^p},\quad p\ge 1, f\in H^1, g\in
H^p.\end{equation*} Actually, (\ref{convdefint})
 implies that if $f\in H^1$ and $g\in \hol (\mathbb D)$, then
\begin{equation}\label{convh1hp} M_p(r, f\star g)\,\le \Vert f\Vert _{H^1}M_p(r, g), \quad 0<r<1.\end{equation}
Also, taking $\rho =r$ in (\ref{convdefint}), we obtain
\begin{equation}\label{Mprsquare}M_p(r^2, f\star g)\,\le M_1(r,
f)M_p(r, g),\quad f, g\in \hol (\mathbb D),\,\,0<r<1.\end{equation}
\par We shall use also the following elementary result.
\begin{lemma}\label{convHpAp} Suppose that $\alpha >-1$, $1\le
p<\infty $, $f\in H^p$, and $g\in A^p_\alpha $. Then
$$\Vert f\star g\Vert _{A^p_\alpha }\,\lesssim \Vert f\Vert
_{H^1}\Vert g\Vert_{A^p_\alpha }.$$
\end{lemma}
\begin{pf}
\begin{align*}\Vert f\ast g\Vert _{A^p_\alpha }^p &\asymp
\int_0^1(1-r)^\alpha M_p^p(r, f\ast g)\,dr\,\lesssim
\int_0^1(1-r)^\alpha \Vert f\Vert _{H^1}^pM_p^p(r,g)\,dr\asymp \Vert
f\Vert_{H^1}^p\Vert g\Vert _{A^p_\alpha }^p.
\end{align*}
\end{pf}
\par\bigskip
\section{The case $p=2$}\label{pequal2}
The results stated for $p=2$ are essentially contained in
\cite{GGPr} and \cite{Bl-Ga-Gi}.
\par Indeed, Theorem\,\@1 of \cite{GGPr} and
Theorem\,\@\ref{th:mlip} give that $\mathcal R_{(\eta )}\in \mathcal
B (H^2)$ if and only if $F_{(\eta )}\in \Lambda^2_{1/2}$.
\par Similarly, part\,\@(d) of Theorem\,\@\ref{compactnessHp}
follows
from Theorem\,\@2 of \cite{GGPr} and the \lq\lq little oh\rq\rq
analogue of Theorem\,\@\ref{th:mlip}.
\par\medskip Notice that, for $f\in \hol (\mathbb D)$,
$f(z)=\sum_{n=0}^\infty a_nz^n$ ($z\in \mathbb D$), we have that
$$\Vert f\Vert _{\mathcal D^2_\alpha }^2\,\asymp\,\vert a_0\vert
^2\,+\,\sum_{n=1}^\infty n^{1-\alpha }\vert a_n\vert ^2.$$ Then
Theorem\,\@3 of \cite{Bl-Ga-Gi} shows that, for $\alpha >0$,
\begin{equation}\label{condD2alpha}\mathcal R_{(\eta )}\in \mathcal B(\mathcal D^2_\alpha
)\,\,\Leftrightarrow\,\,\sum_{n=N}^\infty n^{1-\alpha }\vert
\eta_n\vert ^2\,=\,\og (N^{-\alpha }),\,\,\,\text{as $N\to\infty
$}.\end{equation} Now,
$$\sum_{n=1}^Nn^2\vert \eta_n\vert^2\,=\,\sum_{n=1}^Nn^{1+\alpha }n^{1-\alpha }\vert
\eta_n\vert^2 $$ and then, using the Lemma in p.\,\@101 of
\cite{Du:Hp}, we see that
$$\sum_{n=N}^\infty n^{1-\alpha }\vert
\eta_n\vert ^2\,=\,\og (N^{-\alpha }),\,\,\text{as $N\to\infty
$}\,\,\,\Leftrightarrow\,\,\,\sum_{n=1}^Nn^2\vert
\eta_n\vert^2\,=\,\og (N),\,\,\text{as $N\to\infty $}.$$ This,
(\ref{condD2alpha}), and Theorem\,\@\ref{th:mlip} prove part (v) of
Theorem\,\@\ref{Dpa}. The \lq\lq little oh\rq\rq \, analogues of the
results we have just used give part (vi).
\par\smallskip
Parts (v) and (vi) of Theorem\,\@\ref{Apa} follow from parts (v) and
(vi) of Theorem\,\@\ref{Dpa} and the fact that $A^2_\alpha =\mathcal
D^2_{\alpha +2}$ (see Lemma\,\@3.\,\@10 of \cite{Zhu}).
\par\medskip \begin{remark} For $-1<\alpha \le 0$,  the sequences $(\eta )$ for which
$\mathcal R_{(\eta )}$ is bounded (compact) on $\mathcal D^2_{\alpha
}$ are characterized in Theorems 1 and 2 of
\cite{Bl-Ga-Gi}.\end{remark}
\par\bigskip
\section{Necessary conditions}\label{necessary}

\par\medskip
In this section we are going to prove part (a) of theorems
\ref{boundedness-Hp} and \ref{compactnessHp} and parts (iii) and
(iv) of theorems \ref{Apa} and \ref{Dpa}. \par \medskip We are going
to use the following result. \begin{proposition}\label{varphiN}For
$N\,=\,2, 3, 4, \dots $, set $a_N\,=\,1-\frac{1}{N}$ and
\begin{equation}\label{phiN}\varphi _N(z)\,=\,
\sum_{n=N}^{2N-1}\frac{1} {\sum_{k=1}^nka_N^k}z^n\,\quad z\in
\mathbb D.\end{equation} Then
\begin{equation}\label{normvarphiN}\Vert\varphi_N\Vert
_{H^1}\lesssim \frac{\log N}{N^2}.\end{equation}
\end{proposition}
\begin{pf}
Notice that
$$\varphi_N(z)\,=\,\sum_{n=N}^{2N-1}\frac{1}
{\sum_{k=1}^nka_N^k}z^n\,=\,z^N\psi_N(z),\quad z\in \mathbb D,$$
where
$$ \psi_N(z)\,=\,\sum_{n=0}^{N-1}\frac{1}{\sum_{k=1}^{n+N}ka_N^k}z^n,\quad z\in \mathbb
D.$$ Now, the sequence of Taylor coefficients of $\psi _N$ is
decreasing. Then using \cite[Theorem\,\@1.\,\@1]{Pav-dec}, we obtain

$$\Vert
\varphi_N\Vert_{H^1}\,=\,\Vert
\psi_N\Vert_{H^1}\,\asymp\,\sum_{n=0}^{N-1}\frac{1}{(n+1)\sum_{k=1}^{n+N}ka_N^k}\,\le\,
\sum_{n=0}^{N-1}\frac{1}{(n+1)\sum_{k=1}^{N}ka_N^k}\,\asymp\,\frac{\log
N}{N^2}.$$
\end{pf}

\par For $1\le p<\infty $ and $N\,=\,2, 3, 4,\dots $
set
\begin{equation}\label{deffN}f_{N,p}(z)\,=\,\frac{1}{N^{2-\frac{1}{p}}}\frac{a_Nz}{(1-a_Nz)^2}\,
=\,\frac{1}{N^{2-\frac{1}{p}}}\sum_{n=1}^\infty na_N^nz^n,\quad z\in
\mathbb D.\end{equation}
\begin{Pf}{\it Theorem\,\@\ref{boundedness-Hp}\,\@(a) and Theorem\,\@\ref{compactnessHp}\,\@(a).}
Assume that $1<p<\infty $ and $\mathcal R_{(\eta )}\in \mathcal
B(H^p)$. As observed in Section\,\@\ref{intro} this implies that
$F_{(\eta )}$ is a well defined analytic function in $\mathbb D$. In
fact, since $F_{(\eta )}=\mathcal R_{(\eta )}(1)$, $F_{(\eta )}\in
H^p$.
\par We have that $f_{N,p}\in
H^p$ for all $N$ and
\begin{equation}\label{normfN}\Vert f_{N,p}\Vert _{H^p}\lesssim
1.\end{equation} Since $\mathcal R_{(\eta )}\in \mathcal B(H^p)$,
\begin{equation}\label{normCfN}\Vert \mathcal R_{(\eta )}(f_{N,p})\Vert _{H^p}\lesssim \Vert f_{N,p}\Vert
_{H^p}.\end{equation} Now,
$$
\mathcal R_{(\eta
)}(f_{N,p})(z)\,=\,\frac{1}{N^{2-\frac{1}{p}}}\sum_{n=1}^\infty
\eta_n\left (\sum_{k=1}^nka_N^k\right )z^n,\quad z\in \mathbb D.$$
Then it follows that
\begin{equation}\label{DeltaNF}\Delta_N(F_{(\eta
)})(z)\,=\,\sum_{n=N}^{2N-1}\eta_nz^n\,=\,N^{2-\frac{1}{p}}\mathcal
R_{(\eta )}(f_{N,p})\star \varphi _N(z),\quad z\in \mathbb
D,\,\,N\ge 2.\end{equation} Using Proposition\,\@\ref{varphiN}, we
obtain
\begin{align}\label{normANF}\Vert \Delta_N(F_{(\eta
)})\Vert _{H^p}\,&\lesssim \,N^{2-\frac{1}{p}}\Vert R_{(\eta
)}(f_{N,p})\Vert_{H^p}\Vert \varphi_N\Vert _{H^1}\,\lesssim
\,N^{-\frac{1}{p}}\log N\Vert R_{(\eta
)}(f_{N,p})\Vert_{H^p}.\end{align} Then, using Theorem\,\@2 of
\cite{GiGo}, (\ref{lampomega}), and (\ref{normCfN}), this yields
that $F_{(\eta )}$ satisfies (\ref{condLomega}).
\par\smallskip
Suppose now that $\mathcal R_{(\eta )}\in \mathcal K(H^p)$. Since
$f_{N,p}\to 0$, as $N\to \infty $, uniformly in compact subsets of
$\mathbb D$, using Lemma\,\@3.\,\@7 of \cite{Tj}, we see that
$\mathcal R_{(\eta )}(f_{N,p})\to 0$, as $N\to \infty $, in $H^p$.
This and  (\ref{normANF}) give that
$$\Vert \Delta _N(F_{(\eta
)})\Vert_{H^p}\,=\,\op \left (\left (N^{-\frac{1}{p}}\right )\log
N\right ),\quad \text {as $N\to \infty $}. $$ This and the \lq\lq
little oh\rq\rq \, version of Theorem\,\@2 of \cite{GiGo} give that
$F_{(\eta )}$ satisfies (\ref{condlomega})
\end{Pf}
\par\medskip
\begin{Pf}{\it Theorem\,\@\ref{Apa}\,\@(iii) and Theorem\,\@\ref{Apa}\,\@(iv).}
Assume that $1<p<\infty $, $-1<\alpha <2p-2$, and $\mathcal R_{(\eta
)}\in \mathcal B(A^p_\alpha )$. \par For $N=2, 3, 4, \dots $, set
$$g_{N,p}(z)\,=\,\frac{1}{N^{2-\frac{\alpha
+2}{p}}}\frac{a_Nz}{(1-a_Nz)^2}\,=\,N^{\frac{\alpha
+1}{p}}f_{N,p}(z),\quad z\in \mathbb D,$$ where $f_{N,p}$ is defined
in (\ref{deffN}). We have that $g_{N,p}\in A^p_\alpha $ for all $N$
and, since $\alpha <2p-2$, Lemma\,\@3.\,\@10 of \cite{Zhu} shows
that and $\Vert g_{N,p}\Vert _{A^p_\alpha }\asymp 1$. Using
(\ref{DeltaNF}), we see that
$$\Delta_N(F_{(\eta )})\,=\,N^{\frac{-(\alpha
+1)}{p}}N^{2-\frac{1}{p}}\mathcal R_{(\eta )}(g_{N.p})\star \varphi
_N.$$ This, Lemma\,\@\ref{convHpAp}, and
Proposition\,\@\ref{varphiN} give
\begin{equation}\label{DeltaNFAp}\Vert \Delta_N(F_{(\eta )})\Vert_{A^p_\alpha
}\,\lesssim N^{\frac{-(\alpha +1)}{p}}N^{2-\frac{1}{p}}\Vert
\mathcal R_{(\eta )}(g_{N.p})\Vert _{A^p_\alpha }\frac{\log
N}{N^2}\,=\, N^{\frac{-(\alpha +1)}{p}}N^{-\frac{1}{p}}\log N\Vert
\mathcal R_{(\eta )}(g_{N.p})\Vert _{A^p_\alpha }.\end{equation}
Using Lemma\,\@3.\,\@1 of \cite{MP} it follows easily that $\Vert
\Delta _N\left (F_{(\eta )}\right
)\Vert_{H^p}\asymp\,N^{\frac{1+\alpha }{p}}\Vert \Delta _N\left
(F_{(\eta )}\right )\Vert_{A^p_\alpha }$. Then it follows that
$\Vert \Delta _N\left (F_{(\eta )}\right )\Vert_{H^p}\lesssim
N^{-\frac{1}{p}}\log N$, and, hence, $F_{(\eta )}$ satisfies
(\ref{condLomega}).
\par\smallskip
Suppose now that $\mathcal R_{(\eta )}\in \mathcal K(A^p_\alpha )$.
Since  $g_{N,p}\,\rightarrow\, 0$, as $N\to\infty $, uniformly in
compact subsets of $\mathbb D$, it follows that $\mathcal R_{(\eta
)}(g_{N,p})\,\rightarrow\, 0$, as $N\to\infty $, in $A^p_\alpha $.
This leads to $\Vert \Delta _N\left (F_{(\eta )}\right
)\Vert_{A^p_\alpha }\,=\op \left (N^{\frac{-(\alpha
+1)}{p}}N^{-\frac{1}{p}}\log N\right )$ and to $\Vert \Delta _N\left
(F_{(\eta )}\right )\Vert_{H^p}=\,\op\left (N^{-\frac{1}{p}}\log
N\right )$. Hence, $F_{(\eta )}$ satisfies (\ref{condlomega}).
\end{Pf}
\par\medskip
\begin{Pf}{\it Theorem\,\@\ref{Dpa}\,\@(iii) and Theorem\,\@\ref{Dpa}\,\@(iv).}
Assume that $1<p<\infty $, $-1<\alpha <2p-2$, and $\mathcal R_{(\eta
)}\in \mathcal B(\mathcal D^p_\alpha )$. \par For $N=2, 3, 4, \dots
$, set
$$h_{N,p}(z)\,=\,\frac{1}{
N^{2-\frac{\alpha +2}{p}}}\frac{a_Nz}{1-a_Nz}\,=\,\frac{1}{
N^{2-\frac{\alpha +2}{p}}}\sum_{n=1}^\infty a_N^nz^n,\quad z\in
\mathbb D.$$ Then $h_{N,p}\in \mathcal D^p_\alpha $ for all $N$ and,
using the fact that $2p-\alpha -2>0$ and Lemma\,\@3.\,\@10 of
\cite{Zhu}, we see that
$$\Vert h_{N,p}\Vert _{\mathcal D^p_\alpha
}^p\,\asymp\,\frac{1}{N^{2p-\alpha -2}}\int_{\mathbb
D}\frac{(1-\vert z\vert )^\alpha }{\vert
1-a_Nz\vert^{2p}}\,dA(z)\,\asymp\,1.$$ This implies that $\Vert
\mathcal R_{(\eta )}(h_{N,p})\Vert _{\mathcal D^p_\alpha }\lesssim
1$. \par In order to simplify the notation, we shall let $S$ denote
the shift operator on $\hol (\mathbb D)$, that is,
\begin{equation}\label{shift}S[f](z)\,=\,zf(z),\quad z\in \mathbb
D,\,\,\,f\in \hol (\mathbb D).\end{equation} We have
$$\mathcal R_{(\eta )}(h_{N,p})(z)\,=\,\frac{1}{N^{2-\frac{\alpha +2}{p}}}\sum_{n=1}^\infty \eta_n\left
(\sum_{k=1}^na_N^k\right )z^n,\quad z\in \mathbb D,$$ and, hence,
$$S\left (R_{(\eta )}(h_{N,p})^\prime \right )(z)\,=\,\frac{1}{N^{2-\frac{\alpha +2}{p}}}\sum_{n=1}^\infty
n\eta _n\left (\sum_{k=1}^na_N^k\right )z^n,\quad z\in \mathbb D.$$
Then it follows that \begin{equation}\label{DeltaNFprime}\Delta
_N\left (S(F_\eta ^\prime )\right )\,=\,N^{2-\frac{\alpha
+2}{p}}S\left (R_{(\eta )}(h_{N,p})^\prime \right )\star
\chi_N,\end{equation} where,
$$\chi_N(z)\,=\,\sum_{n=N}^{2N-1}\frac{1}{\sum_{k=1}^na_N^k}z^n,\quad
z\in \mathbb D.$$ Arguing as in the proof of
Proposition\,\@\ref{phiN}, we deduce that
\begin{equation}\label{normchi}\Vert \chi_n\Vert
_{H^1}\,\asymp\,\frac{\log N}{N}.\end{equation} Using this,
(\ref{DeltaNFprime}), and Lemma\,\@\ref{convHpAp}, we obtain
$$\Vert \Delta_N(F_{(\eta )}^\prime )\Vert _{A^p_\alpha }\,\lesssim
N^{1-\frac{\alpha +2}{p}}\log N\Vert \mathcal R_{(\eta
)}(h_{N,p})\Vert _{\mathcal D^p_\alpha }.$$ Then Lemma\,\@3.\,\@1 of
\cite{MP} implies that $$\Vert \Delta_N(F_{(\eta )}^\prime )\Vert
_{H^p}\,\lesssim N^{1-\frac{1}{p}}\log N\Vert \mathcal R_{(\eta
)}(h_{N,p})\Vert _{\mathcal D^p_\alpha }\,\lesssim
N^{1-\frac{1}{p}}\log N$$ and this gives that $F_{(\eta )}$
satisfies (\ref{condLomega}).
\par\smallskip
Suppose now that $\mathcal R_{(\eta )}\in \mathcal K(\mathcal
D^p_\alpha )$. Since  $h_{N,p}\,\rightarrow\, 0$, as $N\to\infty $,
uniformly in compact subsets of $\mathbb D$, it follows that
$\mathcal R_{(\eta )}(h_{N,p})\,\rightarrow\, 0$, as $N\to\infty $,
in $\mathcal D^p_\alpha $. This leads to $$\Vert \Delta_N(F_{(\eta
)}^\prime )\Vert _{A^p_\alpha }\,=\,\op \left (N^{1-\frac{\alpha
+2}{p}}\log N\right )$$ and to
$$\Vert \Delta_N(F_{(\eta )}^\prime )\Vert
_{H^p}\,=\,\op\left (N^{1-\frac{1}{p}}\log N\right ),$$ which
implies that $F_{(\eta )}$ satisfies (\ref{condlomega}).
\end{Pf}
\par\medskip
\section{Sufficient conditions}\label{sufficient}

In order to prove the sufficiency statements in
Theorem~\ref{boundedness-Hp} and Theorem~\ref{compactnessHp} we will
work with the spaces $\mathcal D^p_{p-1}$ ($1\le p<\infty $) and
with the spaces, which we will call $\mathcal X_{q,p}$, consisting
of those $f\in \hol (\mathbb D)$ whose derivative belongs to the
mixed norm space $H^{q,p,\alpha }$ with $\alpha
=1+\frac{1}{p}-\frac{1}{q}$ ($1\le q\le p<\infty $). Thus,
$$\mathcal X_{q,p}=\left \{ f\in \hol (\mathbb D): \Vert f\Vert _{\mathcal X_{q,p}}\ig
\left (\vert f(0)\vert
^p\,+\,\int_0^1(1-r)^{p(1-\frac{1}{q})}M_q^p(r, f^\prime )\,dr\right
)^{1/p}<\infty \right \}.$$
\par Notice that
\begin{equation}\label{DpXPP}\mathcal X_{p,p}\,=\,\mathcal D^p_{p-1},\quad
1\le p<\infty .\end{equation}
\par It is well known
(see, e.\,\@g., \cite{LP,Lu,Vi}) that
\begin{equation}\label{dppless2}\mathcal D^p_{p-1}\subset H^p,\quad
0<p\le 2,\end{equation} and
\begin{equation}\label{dppg2}H^p\subset \mathcal D^p_{p-1},\quad
2<p<\infty ,\end{equation} the inclusions being continuous. \par  A
result of Flett \cite[Theorem\,\@1]{Flett} (see also
Corollary\,\@3.\,\@1 of \cite{MP}) asserts that
\begin{equation}\label{Xpqemb}\text{If $1<q<p$ then $\mathcal X_{q,p}$ is
continuously embedded in $H^p$.}\end{equation}
\par
Bearing in mind (\ref{dppless2}) and (\ref{Xpqemb}), the sufficiency
statements in Theorem~\ref{boundedness-Hp} and
Theorem~\ref{compactnessHp} will follow readily from the following
result.
\begin{theorem}\label{DpXqp}
Let $(\eta )\,=\,\{ \eta_n\} _{n=0}^\infty $ be a sequence of
complex numbers.
\par (i) If\, $1\le p<\infty $ and $F_{(\eta )}\in \Lambda ^p_{1/p}$,
then $\mathcal R_{(\eta )}$ is a bounded operator from $H^p$ into
$\mathcal D^p_{p-1}$.
\par (ii) If\, $1<p<\infty $ and $F_{(\eta )}\in \lambda ^p_{1/p}$,
then $\mathcal R_{(\eta )}$ is a compact operator from $H^p$ into
$\mathcal D^p_{p-1}$.
\par (iii) If\, $1<q<p<\infty $ and $F_{(\eta )}\in \lambda ^q_{1/q}$,
then $\mathcal R_{(\eta )}$ is a compact operator from $H^p$ into
$\mathcal X_{q,p}$.
\end{theorem}
\par Our proof of this result is inspired in the work of M.~Nowak
\cite{Nowak} who proved that the Ces\`{a}ro operator is bounded from
$H^p$ into $\mathcal D^p_{p-1}$ ($0<p<\infty $). \par Let us notice
that, for $h\in \hol (\mathbb D)$, we have
\begin{align}\label{Xpqrsquare}
\Vert h\Vert _{\mathcal X_{q,p}^p}\,\asymp\,&\vert h(0)\vert
^p\,+\,\int_0^1(1-r)^{p\left (1-\frac{1}{q}\right )}M_q^p(r,
h^\prime )\,dr\\ =\,& \vert h(0)\vert
^p\,+\,\int_0^12r(1-r^2)^{p\left (1-\frac{1}{q}\right )}M_q^p(r^2,
h^\prime )\,dr \nonumber\\ \lesssim\,& \vert h(0)\vert
^p\,+\,\int_0^1(1-r)^{p\left (1-\frac{1}{q}\right )}M_q^p(r^2,
h^\prime )\,dr.\nonumber\end{align}
\par\medskip
\begin{Pf}{\it Theorem\,\@\ref{DpXqp}.} Suppose that $1\le q\le
p<\infty $ and $F_{(\eta )}\in \Lambda ^q_{1/q}$.
\par Take $f\in H^p$, $f(z)=\sum_{n=0}^\infty a_nz^n$ ($z\in \mathbb
D$) and set
\begin{equation}\label{g}g(z)\,=\,\frac{f(z)}{1-z}\,=\,\sum_{n=0}^\infty
\left (\sum_{k=0}^na_k\right )z^n,\quad z\in \mathbb
D.\end{equation} With $S$ denoting the shift operator as above, we
have,
\begin{align*}S\left [\mathcal R_{(\eta )}(f)\right ](z)
\,=\,\sum_{n=0}^\infty (n+1)\eta_n\left
(\frac{1}{n+1}\sum_{k=0}^na_k\right )z^{n+1},\quad z\in \mathbb
D,\end{align*} and then
\begin{align}\label{SRf}S\left [\mathcal R_{(\eta )}(f)\right ]^\prime(z)
\,=\,\sum_{n=0}^\infty (n+1)\eta_n\left (\sum_{k=0}^na_k\right
)z^n\,=\,S\left [F_{(\eta )}\right ]^\prime \star g(z),\quad z\in
\mathbb D.\end{align} Using (\ref{Mprsquare}) and the continuity of
the shift operator, we obtain, for $1/2<r<1$,
\begin{align}\label{MqReta}M_q\left (r^2, \left [\mathcal R_{(\eta )}(f)\right ]^\prime\right
)&\,\lesssim M_q\left (r, S\left [\mathcal R_{(\eta )}(f)\right ]
^\prime \right )
\\ &\,\lesssim M_q\left (r, S[F_{(\eta )}]^\prime \right
)M_1(r, g)\nonumber\\ &\,\lesssim M_q\left (r, F_{\eta }^\prime
\right )M_1(r, g)\nonumber\\ \,& \lesssim (1-r)^{\frac{1}{q}-1}\beta
_{q,\frac{1}{q}}\left (r, F_{(\eta )}\right )M_1(r, g).\nonumber
\end{align}
Take $\gamma >1$ and let $\sigma $ be the exponent conjugate to
$\gamma p$, that is, $\frac{1}{\gamma p}+\frac{1}{\sigma }=1$. Using
H\"{o}lder`s inequality, we obtain
$$M_1(r, g)=\frac{1}{2\pi }\int_0^{2\pi }\left \vert
\frac{f(re^{i\theta }}{1-re^{i\theta }}\right \vert d\theta  \,\le
M_{\gamma p}(r, f)\left (\frac{1}{2\pi }\int_0^{2\pi } \frac{d\theta
}{\vert 1-re^{i\theta } \vert ^\sigma }\right )^{1/\sigma } \lesssim
\frac{1}{(1-r)^{\frac{1}{\gamma p}}}M_{\gamma p}(r, f).
$$
This and (\ref{MqReta}) give
\begin{align*}M_q\left (r^2, \left [\mathcal R_{(\eta )}(f)\right ]^\prime\right
)&\,\lesssim (1-r)^{\frac{1}{q}-1}\beta _{q,\frac{1}{q}}\left (r,
F_{(\eta )}\right )\frac{1}{(1-r)^{\frac{1}{\gamma p}}}M_{\gamma
p}(r, f).
\end{align*}
This and (\ref{Xpqrsquare})  imply that

\begin{align}\label{norXpqReta}\Vert \mathcal R_{(\eta)}(f)\Vert _{\mathcal X_{q,p}}^p\,\lesssim &
\vert \eta_0\vert ^p\vert a_0\vert ^p\,+\int_0^1(1-r)^{p\left (1-\frac{1}{q}\right )}M_q^p
\left (r^2, \left [\mathcal R_{(\eta )}(f)\right ]^\prime\right )\,dr\\
\lesssim &\vert \eta_0\vert ^p\Vert f\Vert
_{H^p}^p\,+\,\int_0^1\beta ^p _{q,\frac{1}{q}}\left (r, F_{(\eta
)}\right )\frac{1}{(1-r)^{\frac{1}{\gamma }}}M_{\gamma p}^p(r,
f)dr.\nonumber
\end{align}
Since $F_{(\eta )}\in \Lambda ^q_{1/q}$,\,\, $\sup_{0<r<1}\beta _{q,
\frac{1}{q}}\left (r, F_{(\eta )}\right )<\infty $. Then, using
Theorem\,\@5.\,\@11 of \cite{Du:Hp}, it follows that
\begin{equation}\label{boubdeXqp}\Vert \mathcal R_{(\eta)}(f)\Vert _{\mathcal
X_{q,p}}^p\lesssim \Vert f\Vert _{H^p}^p.\end{equation} Thus we have
proved that $\mathcal R_{(\eta )}$ is bounded from $H^p$ into
$\mathcal X_{q,p}$.
\par\medskip Part (i) of Theorem\,\@\ref{DpXqp} follows taking $q=p$ and
using (\ref{DpXPP}).
\par\smallskip Parts (ii) and (iii) of  Theorem\,\@\ref{DpXqp} follow from the
following result.
\begin{proposition}\label{compact} Let $(\eta )=\{ \eta_n\}
_{n=0}^\infty $ be a sequence of complex numbers. If $1<q\le
p<\infty $ and $F_{(\eta )}\in \lambda ^q_{1/q}$, then $\mathcal
R_{(\eta )}$ is compact form $H^p$ into $\mathcal X_{q,p}$.
\end{proposition}
\begin{pf} Suppose that $1<q\le p<\infty
$ and $F_{(\eta )}\in \lambda ^p_{1/p}$. We already know that
$\mathcal R_{(\eta )}$ is bounded from $H^p$ into $\mathcal
X_{q,p}$. For $N=2, 3, 4, \dots $ let $\left (\eta ^{(N)}\right
)\,=\,\{ \eta^{(N)}_n\} _{n=0}^\infty $ be the sequence defined by
$\eta ^{(N)}_n=0$, if $0\le n\le N$, and $\eta ^{(N)}_n=\eta_n$, if
$n\ge N+1$. Set also $\mathcal R_N:\,\hol (\mathbb D)\rightarrow
\hol (\mathbb D)$ be the operator defined by
\begin{equation}\label{defRN}\mathcal R_N(f)(z)\,=\,\sum_{n=0}^N\eta _n\left (\sum_{k=0}^na_k\right
)z^n, \quad z\in \mathbb D,\,\,\,f\in \hol (\mathbb
D).\end{equation} The operators $\mathcal R_N$ are finite rank
operators from $H^p$ into $\mathcal X_{q,p}$. Then it suffices to
prove that
\begin{equation}\label{oper-norm}\mathcal R_N\,\rightarrow\, \mathcal R_{(\eta
)},\quad\text{as $N\to \infty $,\, in the operator
norm}.\end{equation} Notice that
$$\mathcal R_{(\eta )}\,-\,\mathcal R_N\,=\,\mathcal R_{(\eta
^{(N)})}\,\,\,\text{and}\,\,\,F_{(\eta ^{(N)})}\,=\,F_{(\eta
)}\,-\,S_NF_{(\eta )}.$$ Take $\gamma >1$. Then, using
(\ref{norXpqReta}) with $(\eta ^{(N)})$ in the place of $(\eta )$,
and Theorem\,\@5.\,\@11 of \cite{Du:Hp} as above, we see that if
$f\in H^p$ then \begin{align*}\Vert \left (\mathcal R_{(\eta
)}\,-\,\mathcal R_N\right )(f)\Vert _{\mathcal X_{q,p}}^p\,\lesssim
\,&\int_0^1\beta _{q,\frac{1}{q}}^p\left (r, F_{(\eta
)}\,-\,S_NF_{(\eta )}\right
)\frac{1}{(1-r)^{\frac{1}{\gamma }}}M_{\gamma p}^p(r, f)dr\\
\,\lesssim &\,\beta_{q,\frac{1}{q}}^q\left (F_{(\eta )}-S_nF_{(\eta
)}\right )\int_0^1\frac{1}{(1-r)^{\frac{1}{\gamma }}}M_{\gamma
p}^p(r, f)dr\\ \,\lesssim &\,\beta_{q,\frac{1}{q}}^q\left (F_{(\eta
)}-S_nF_{(\eta )}\right )\Vert f\Vert _{H^p}^p.\end{align*} Using
Proposition\,\@\ref{partialsumlamba}, (\ref{oper-norm}) follows as
desired.
\end{pf}
\end{Pf}
\par\medskip

Since, for $\alpha >-1$ and $p>0$, $A^p_\alpha =\mathcal D^p_{\alpha
+p}$ (see Theorem\,\@4.\,\@28 of \cite{Zhu}), parts (i) and (ii) of
Theorem\,\@\ref{Apa} follow from parts (i) and (ii) of
Theorem\,\@\ref{Dpa}.

\par\medskip

\begin{Pf}{\it Theorem\,\@\ref{Dpa}\,\@(i) and Theorem\,\@\ref{Dpa}\,\@(ii).}
Assume that  $1\le p<\infty $, $\alpha >p-2$, and $F_{(\eta )}\in
\Lambda ^p_{1/p}$. Using Theorem\,\@4.\,\@28 of \cite{Zhu}, we have
to prove that, for $f\in \mathcal D^p_\alpha $,
\begin{equation}\label{intdobleprima}\int_{\mathbb D}(1-\vert z\vert
)^{\alpha +p}\vert \mathcal R_{(\eta )}(f)^{\prime\prime }(z)\vert
^p\,dA(z)\,\lesssim \Vert f\Vert _{\mathcal D^p_\alpha
}^p.\end{equation} This is equivalent to
\begin{equation}\label{intrdobleprima}\int_0^1(1-r
)^{\alpha +p}M_p^p(r^4, \mathcal R(f)^{\prime\prime })\,dr\,\lesssim
\Vert f\Vert _{\mathcal D^p_\alpha }^p.\end{equation}

So, take $f\in \mathcal D^p_\alpha $, $f(z)=\sum_{n=0}^\infty
a_nz^n$ ($z\in \mathbb D$). Just as in the proof of
Theorem\,\@\ref{DpXqp}, set $g(z)=\frac{f(z)}{1-z}$ ($z\in \mathbb
D$). Bearing in mind (\ref{SRf}) and the continuity of the shift
operator on $H^p$, we see that (\ref{intrdobleprima}) will follow
from
\begin{equation}\label{Fprimegprime}\int_0^1(1-r
)^{\alpha +p}M_p^p(r^4, F_{(\eta )}^\prime \star g^\prime
)\,dr\,\lesssim \Vert f\Vert _{\mathcal D^p_\alpha
}^p.\end{equation} Notice that
$$g^\prime (z)\,=\,\varphi (z)\,+\,\phi (z),\quad z\in \mathbb
D,$$ with
$$\varphi (z)\,=\,\frac{f^\prime (z)}{1-z},\quad \phi
(z)\,=\,\frac{f(z)}{(1-z)^2},\quad z\in \mathbb D.$$ Using this,
(\ref{Mprsquare}), and the fact that $F_{(\eta )}\in
\Lambda^p_{1/p}$, we see that
\begin{align}\label{r4Fprimegprime}&\int_0^1(1-r
)^{\alpha +p}M_p^p(r^4, F_{(\eta )}^\prime \star g^\prime
)\,dr\lesssim \int_0^1(1-r)^{\alpha +p}M_p^p(r^2, F_{(\eta )}^\prime
)M_1^p(r^2, g^\prime )\,dr\\&\lesssim
\beta_{p,\frac{1}{p}}^p(F_{(\eta )})\left [\int_0^1(1-r)^{\alpha
+1}M_1^p(r^2, \varphi)\,dr+\int_0^1(1-r)^{\alpha +1}M_1^p(r^2,
\phi)\,dr\right ]\nonumber
\\&\lesssim \beta_{p,\frac{1}{p}}^p(F_{(\eta )})\left
[ \int_0^1(1-r)^{\alpha +1}M_1^p(r,
\varphi)\,dr+\int_0^1(1-r)^{\alpha +1}M_1^p(r^2, \phi)\,dr\right
].\nonumber
\end{align}
Using H\"{o}lder's inequality, we see that
$$M_1(r, \varphi )=\frac{1}{2\pi }\int_0^{2\pi }\frac{\vert
f^\prime (re^{i\theta })\vert }{\vert 1-re^{i\theta }\vert
}\,d\theta \lesssim \frac{1}{(1-r)^{1/p}}M_p(r, f^\prime ).$$ This
implies that
\begin{align}\label{r4varphi}\int_0^1(1-r)^{\alpha +1}M_1^p(r,
\varphi)\,dr& \,\lesssim \int_0^1(1-r)^\alpha M_p^p(r, f^\prime
)\,dr\,\lesssim \Vert f\Vert _{\mathcal D^p_{\alpha }}^p.\end{align}
Take now $\gamma >1$ such that $\alpha >p-2+\frac{1}{\gamma }$. Let
$\sigma $ be the exponent conjugate to $\gamma p$, that is,
$\frac{1}{\gamma p}+\frac{1}{\sigma }=1$. Using H\"{o}lder's
inequality again, we see that
$$M_1(r^2, \phi )\,\lesssim M_{\gamma p}(r^2, f)\left (\int_0^{2\pi
}\frac{d\theta }{\vert 1-re^{i\theta }\vert ^{2\sigma }}\right
)^{1/\sigma }\asymp\frac{1}{(1-r)^{1+\frac{1}{\gamma p}}}M_{\gamma
p}(r^2, f)$$ and then it follows that
$$\int_0^1(1-r)^{\alpha +1}M_1^p(r^2,
\phi)\,dr \,\lesssim \int_0^1(1-r)^{\alpha +1-p-\frac{1}{\gamma
}}M_{\gamma p}^p(r^2, f)\,dr.$$ Since $\alpha +1-p-\frac{1}{\gamma
}>-1$, using Theorem\,\@5.\,\@6 of \cite{Du:Hp}, we obtain
$$\int_0^1(1-r)^{\alpha +1}M_1^p(r^2,
\phi)\,dr \,\lesssim \int_0^1(1-r)^{\alpha +1-\frac{1}{\gamma
}}M_{\gamma p}^p(r^2, f^\prime )\,dr.$$ Now, a slight modification
of the proof of Theorem\,\@5.\,\@9 of \cite{Du:Hp} gives that
$$M_{\gamma p}(r^2, f^\prime )\lesssim
\frac{1}{(1-r)^{\frac{1}{p}-\frac{1}{\gamma p}}}M_p(r, f^\prime ).$$
Then we obtain
$$\int_0^1(1-r)^{\alpha +1}M_1^p(r^2,
\phi)\,dr \,\lesssim \int_0^1(1-r)^{\alpha }M_p^p(r, f^\prime
)\,dr\,\lesssim \Vert f\Vert _{\mathcal D^p_\alpha }^p.$$ This,
(\ref{r4varphi}), and (\ref{r4Fprimegprime}) give
\begin{equation}\label{condpa}\int_0^1(1-r
)^{\alpha +p}M_p^p(r^4, \mathcal R(f)^{\prime\prime })\,dr\,\lesssim
\beta_{p,\frac{1}{p}}^p(F_{(\eta )})\Vert f\Vert _{\mathcal
D^p_\alpha }^p.\end{equation} Hence (\ref{Fprimegprime}) holds as
desired.
\par\smallskip
Suppose now that $1<p<\infty $ and $F_{(\eta )}\in
\lambda^{p}_{1/p}$. For $N=2, 3, 4, \dots $ let  $\mathcal R_N$ be
defined as in (\ref{defRN}). Using (\ref{condpa}) and arguing as in
the proof of Proposition\,\@\ref{compact}, we deduce, for $f\in
\mathcal D^p_\alpha $,
\begin{align*}\left \Vert \left (\mathcal R_{(\eta )}\,-\,\mathcal R_N\right
)(f)\right \Vert _{\mathcal D^p_\alpha }^p \lesssim
\beta^p_{p,\frac{1}{p}}\left (F_{(\eta )}\,-\,S_NF_{(\eta )}\right
)\Vert f\Vert _{\mathcal D^p_\alpha }^p.
\end{align*}
 Using
Proposition\,\@\ref{partialsumlamba} we deduce that $\{ \mathcal
R_{(\eta )}-\mathcal R_N\} \,\rightarrow \, 0$, as $N\to \infty $,
in the norm of $\mathcal B(\mathcal D^p_\alpha )$. Since the
operators $\mathcal R_N:\mathcal D^p_\alpha \rightarrow \mathcal
D^p_\alpha $ are finite rank operators, it follows that $\mathcal
R_{(\eta )}\in \mathcal K(\mathcal D^p_\alpha )$. Thus, (ii) is
proved.

\end{Pf}
\par\bigskip
\section{Rhaly operators acting on $H^1$ and some further results}\label{further}
We start this section obtaining a number of results about the action
of the Rhaly operators on $H^1$. \par  In part (i) of
Theorem\,\@\ref{DpXqp} we proved that if $F_{(\eta )}\in \Lambda
^1_1$ then $\mathcal R_{(\eta )}\in \mathcal B(H^1)$. In our next
theorem we improve this result, showing that the condition $F_{(\eta
)}\in \Lambda ^1_1$ actually implies that $\mathcal R_{(\eta )}\in
\mathcal K(H^1)$. We give also two conditions which are necessary
for the boundedness of the Rhaly operator $\mathcal R_{(\eta )}$ in
$H^1$. Before stating the theorem, let us notice that $$\Lambda
^1_1\,=\,\{ f\in \hol (\mathbb D): f^\prime \in H^1\}\,=\,S^1.$$
Also
$$\beta_{1,1}(r, f)\,=\,M_1(r, f^\prime ),\quad
\beta_{1,1}(f)\,=\,\Vert f^\prime\Vert _{H^1}.$$
\begin{theorem}\label{RetaH1} Let $(\eta )=\{ \eta_n\} _{n=0}^\infty
$ be a sequence of complex numbers.
\par (i)\,If $F_{(\eta )}^\prime \in H^1$ then $\mathcal R_{(\eta )}\in \mathcal K(H^1,\mathcal D^1_0)$  and $\mathcal R_{(\eta
)}\in \mathcal K(H^1)$.
\par (ii) If \, $\mathcal R_{(\eta )}\in \mathcal B(H^1)$ then $\sum_{n=0}^Nn\vert \eta_n\vert \,=\,\og (N)$
\par (iii) If\, $\mathcal R_{(\eta )}\in \mathcal B(H^1)$ then $\Vert \Delta _N(F_{(\eta )})^\prime\Vert _{H^1}\,=\,\og (\log N)$,
as $N\to \infty $.
\end{theorem}
\begin{pf} Assume that $F_{(\eta )}^\prime \in H^1$. We already know
that this implies that $\mathcal R_{(\eta )}$ is a bounded operator
from $H^1$ into $\mathcal D^1_0$.
\par It is well known that there are functions $h$ in $H^1$ for which
the sequence of partial sums $S_Nh$ does not converge. However, the
following is true:
\par\smallskip {\it{If $h\in H^1$ then there exists a subsequence $\{
S_{N_k}h\} $ of the sequence of partial sums of $h$ such that
$S_{N_k}h\to h$, as $k\to \infty $, in the norm of $H^1$.}}
\par
For instance, this can be deduced easily from the fact that $\Vert
h_r-h\Vert _{H^1}\to 0$, as $r\to 1$.
\par\smallskip So, take a subsequence $\{ S_{N_k}(F_{(\eta )}^\prime )\}$
of the sequence of the partial sums of $F_{(\eta )}^\prime $
converging to $F_{(\eta )}^\prime $ in $H^1$. Arguing as in the
proof of Proposition\,\@\ref{compact} we obtain that, if $\gamma
>1$ and the operators $\mathcal R_{N_{k+1}}$ are defined as in (\ref{defRN}),
then, for $f\in H^1$, we have \begin{align*}\Vert \left (\mathcal
R_{(\eta )}\,-\,\mathcal R_{N_{k+1}}\right )(f)\Vert_{\mathcal
D^1_0}\,\lesssim \,& \beta_{1,1}\left (F_{(\eta
)}\,-\,S_{N_{k+1}}F_{(\eta )}\right )\Vert f\Vert _{H^1}\\
\,=\,&\Vert F_{(\eta )}^\prime\,-\,\left (S_{N_{k+1}}F_{(\eta
)}\right )^\prime \Vert_{H^1}\Vert f\Vert _{H^1}
\\ \,=\,&\Vert F_{(\eta
)}^\prime\,-\, S_{N_{k}}\left (F_{(\eta )}^\prime\right )
\Vert_{H^1}\Vert f\Vert _{H^1} .\end{align*} Then we see that
$\mathcal R_{(\eta )}\,-\,\mathcal R_{N_{k+1}}\to 0$, as $k\to
\infty $ in the norm of the space $\mathcal B(H^1,\mathcal D^1_0)$.
Since the operators $\mathcal R_{N_{k+1}}:H^1\,\rightarrow\,\mathcal
D^1_0$ are finite rank operators, it follows that $\mathcal R_{(\eta
)}\in \mathcal K(H^1,\mathcal D^1_0)$ and, then also $\mathcal
R_{(\eta )}\in \mathcal K(H^1)$. Thus (i) is proved.
\par\bigskip
Suppose now that $\mathcal R_{(\eta )}\in \mathcal B(H^1)$. For
$N\,=\,2, 3, 4, \dots $ set $a_N\,=\,1-\frac{1}{N}$ and
$$f_N(z)\,=\,\frac{1}{N}\frac{a_Nz}{(1-a_Nz)^2}\,=\,\frac{1}{N}\sum_{n=1}^\infty
na_N^nz^n,\quad z\in\mathbb D.$$ Then, for all $N$, we have that\,
$f_N\in H^1$ and $\Vert f_N\Vert _{H^1}\le 1$. This implies that
\begin{equation}\label{Rfa}\Vert \mathcal R_{(\eta )}(f_N)\Vert
_{H^1}\lesssim 1.\end{equation} We have
\begin{equation}\label{RetafNp1}\mathcal R_{(\eta
)}(f_N)(z)=\frac{1}{N}\sum_{n=1}^\infty \eta_n\left
(\sum_{k=1}^nka_N^k\right )z^n,\quad z\in \mathbb D.\end{equation}
Using Hardy's inequality (see section\,\@3.\,\@6 of \cite{Du:Hp}),
we obtain
$$\frac{1}{N}\sum_{n=1}^N\frac{\vert \eta_n\vert }{n+1}\left
(\sum_{k=1}^nka_N^k\right )\le \pi \Vert \mathcal R_{(\eta
)}(f_N)\Vert _{H^1}\,\lesssim 1.$$ Since there exists  a positive
constant $A$ such that  $a_N^k\ge A$ for $k=1, 2, \dots , N$, we
obtain $\frac{1}{N}\sum_{n=0}^Nn\vert \eta_n\vert \lesssim 1$. Hence
(ii) is proved.
\par\medskip
Using (\ref{RetafNp1}) we see that
$$\Delta_N(F_{(\eta )})\,=\,N\mathcal R_{(\eta )}(f_N)\ast \varphi
_N$$ where $\varphi_N$ is defined in (\ref{phiN}). Using this,
(\ref{Rfa}), and Proposition\,\@\ref{varphiN}, we obtain that $\Vert
\Delta _N(F_{(\eta )})\Vert _{H^1}\,\lesssim \,\frac{\log N}{N}$.
Then, using Lemma\,\@3.\,\@4 of \cite{BSS}, we obtain $$\Vert \Delta
_N(F_{(\eta )})^\prime \Vert _{H^1}\,\lesssim \,\log N.$$ Hence,
(iii) is proved.

\end{pf}
\par\bigskip
Let us turn to prove Theorem\,\@\ref{og1n}.\par\smallskip
\begin{Pf}{\it Theorem\,\@\ref{og1n}.}\par If\, $\eta_n\,=\,\og \left (\frac{1}{n}\right )$, as $n\to \infty
$ then
$$M_2\left (r, F_{(\eta )}^\prime \right )\,=\,\left (\sum_{n=1}^\infty n^2\vert \eta_n\vert
^2r^{2n-2}\right )^{1/2}\,\lesssim \left (\sum_{n=1}^\infty
r^{2n-2}\right )^{1/2}\,\lesssim \frac{1}{(1-r)^{1/2}}.$$ This
implies that $F_\mu \in \Lambda ^2_{1/2}$ and then it follows that
$F_\mu \in \Lambda ^q_{1/q}$ for all $q\in [2, \infty )$. Using
Theorem\,\@\ref{boundedness-Hp}, (i) follows. \par\medskip We shall
use Rademacher functions and Khinchine's inequality to prove (ii).
Let $\{ r_k\}_{k=0}^\infty $ be the sequence of Rademacher
functions, that is,
$$
r_k(t)\,=\,\left\{%
\begin{array}{ll}
   1,\,\,\, & \hbox{if $\frac{2n}{2^{k+1}}\le t<\frac{2n+1}{2^{k+1}}$,\,\,$n=\,0, 1, 2, \dots , 2^k-1$}.\\
   \\
   -1,\,\,\, & \hbox{if $\frac{2n+1}{2^{k+1}}\le t< \frac{2n+2}{2^{k+1}}$,\,\,\,$n=\,0, 1, 2, \dots , 2^k-1$}.
\end{array}%
\right.
$$
Khinchine's inequality is the following (see, e.\,\@g.,
Theorem\,\@8.\,\@4 in chapter V of \cite{Zy} or Appendix A of
\cite{Du:Hp}).
\begin{other}\label{Khin}{\rm[Khinchine's inequality].} For $0<p<\infty $ there
exist constants $A_p$, $B_p$ with $0<A_p<B_p<\infty $ such that for
all natural numbers $m$ and all $m+1$-tuples of complex numbers
$c_0, c_1, \dots , c_m$, we have $$A_p\left (\sum_{j=0}^m\vert
c_j\vert ^2\right )^{p/2}\,\le \int_0^1\left \vert
\sum_{j=0}^mc_jr_j(t)\right \vert ^p\,dt\,\le B_p\left
(\sum_{j=0}^m\vert c_j\vert ^2\right )^{p/2}.$$
\end{other}
\par\medskip Let $P$ is a complex polynomial,
$P(z)=\sum_{j=0}^ma_jz^j$. We define
$$P_t(z)\,=\,\sum_{j=0}^mr_j(t)a_jz^j,\quad 0\le t<1.$$
Using Khinchine's inequality with $c_j=a_je^{ij\theta }$ ($0\le j\le
m$), we obtain
$$A_p\left (\sum_{j=0}^m\vert
a_j\vert ^2\right )^{p/2}\,\le \int_0^1\vert P_t(e^{i\theta })\vert
^p\,dt\,\le \,B_p\left (\sum_{j=0}^m\vert a_j\vert ^2\right
)^{p/2},\quad \theta \in \mathbb R.$$ This implies that $$A_p\left
(\sum_{j=0}^m\vert a_j\vert ^2\right )^{p/2}\,\le \frac{1}{2\pi
}\int_0^{2\pi }\int_0^1\vert P_t(e^{i\theta })\vert ^p\,dt\,d\theta
\,\le \,B_p\left (\sum_{j=0}^m\vert a_j\vert ^2\right )^{p/2}.$$
Interchanging the order of integration in the second term, it
follows that
$$A_p\left
(\sum_{j=0}^m\vert a_j\vert ^2\right )^{p/2}\,\le \int_0^1\Vert
P_t\Vert _{H^p}^p\,dt\,\le \,B_p\left (\sum_{j=0}^m\vert a_j\vert
^2\right )^{p/2}.$$ Clearly, this gives the following:
\begin{equation}\label{texists} \text{There exists $t\in [0, 1)$ such that $\Vert
P_t\Vert _{H^p}\ge A_p\left (\sum_{j=0}^m\vert a_j\vert ^2\right
)^{1/2}.$}\end{equation}
\par\medskip Suppose now that $1\le p<2$.
\par Set $h(z)\,=\,\log\frac{1}{1-z}\,=\,\sum_{n=1}^\infty
\frac{z^n}{n}$ ($z\in \mathbb D$). For all $k$, we have
$$\left (\Delta_{2^k}h\right )^{\prime }(z)\,=\,\sum_{j=2^k}^{2^{k+1}-1}z^{j-1},\quad z\in \mathbb D,\,\,\,k=0, 1, 2,\,\dots .$$

Bearing in mind (\ref{texists}), for each $k$, let us pick $t_k\in
[0, 1)$ such that
$$\Vert \left (\Delta _{2^k}h\right )^\prime _{t_k}\Vert \ge
A_p2^{k/2}.$$ Set
$$\Upsilon (z)\,=\,\sum_{k=0}^\infty \left (\Delta _{2^k}h\right
)_{t_k}(z),\quad z\in \mathbb D.$$ Clearly, $\Upsilon \in \hol
(\mathbb D)$ and $\Delta _{2^k}\Upsilon =\left (\Delta_{2^k}h\right
)_{t_k}$, for all $k$. Let $(\eta )=\{ \eta _n\}_{n=0}^\infty $ be
defined by $\eta_n=\frac{\Upsilon^{(n)}(0)}{n!}$. Then $\Upsilon
=F_{(\eta )}$. We have that $\vert \eta _n\vert =\frac{1}{n}$ for
all $n\ge 1$. Also \begin{equation}\label{nobounded}\Vert \left
(\Delta_{2^k}F_{(\eta )}\right )^\prime \Vert _{H^p}\,=\,\Vert \left
(\Delta_{2^k}h\right )^\prime \Vert _{H^p}\ge
A_p2^{k/2}.\end{equation}
\par If $1<p<2$, then (\ref{nobounded}) implies that $F_{(\eta )}$ does not satisfy (\ref{condLomega}) and then, using Theorem\,\@\ref{boundedness-Hp},
we deduce that $\mathcal R_{(\eta )}\not\in \mathcal B(H^p)$.
\par If $p=1$,  then (\ref{nobounded}) and part (iii) of
Theorem\,\@\ref{RetaH1} show that $\mathcal R_{(\eta )}\not\in
\mathcal B(H^1)$.
\par\medskip Part (iii) of Theorem\,\@\ref{og1n} follows immediately
from the following result (see \cite{GiMer} for $p\ge 2$, and
\cite[Lemma\,\@1]{Mer} for $1<p<\infty $).
\begin{otherp}\label{lamp1n} If\, $1<p<\infty $ and $(\eta )=\{ \eta
_n\}_{n=0}^\infty $ is a decreasing sequence of non negative real
numbers, then $F_{(\eta )}\in \Lambda^p_{1/p}$ if and only if
$\eta_n=\og \left (\frac{1}{n}\right )$, as $n\to \infty $.
\end{otherp}
\end{Pf}
\par\smallskip
For $1\le p\le 2$ we can prove the following.
\begin{proposition}\label{vertetan} Suppose that $1\le p\le 2$ and
let $(\eta )=\{ \eta_n\} _{n=0}^\infty $ be a sequence of complex
numbers such that the sequence $\{ \vert \eta_n\vert \}
_{n=0}^\infty $ is decreasing.
\par If $\mathcal R_{(\eta )}\in \mathcal B(H^p)$ then $\eta_n=\og \left (\frac{1}{n}\right
)$, as $n\to \infty $.
\end{proposition}
\begin{pf}

If $p=1$ the result follows readily from part (ii) of
Theorem\,\@\ref{RetaH1}.
\par
 Assume that $\{ \vert \eta_n\vert \} _{n=0}^\infty $ is
decreasing, $1<p\le 2$, and $\mathcal R_{(\eta )}\in \mathcal
B(H^p)$. For $N=2, 3, 4, \dots $, set $a_N=1-\frac{1}{N}$ and let
$f_N$ be defined as in (\ref{deffN}). We have
$$\mathcal R_{(\eta )}(f_N)(z)
=\frac{1}{N^{2-\frac{1}{p}}}\sum_{n=1}^\infty \eta_n\left
(\sum_{k=0}^nka_N^k\right )z^n,\quad z\in \mathbb D,$$ and $\Vert
R_{(\eta )}(f_N)\Vert _{H^p}\lesssim 1$. Using this, the facts that
$\{ \vert \eta _n\vert \} $ is decreasing and that there exists
$A>0$ such that $a_N^k\ge A$ for $0\le k\le N$, and
Theorem\,\@6.\,\@2 of \cite{Du:Hp}, we deduce that $$\vert
\eta_N\vert ^p\frac{1}{N^{2p-1}}\sum_{n=1}^N(n+1)^{p-2}\left
(\sum_{k=0}^nk\right )^p\lesssim \Vert \mathcal R_{(\eta
)}(f_N)\Vert _{H^p}^p\lesssim 1.$$ Since
$$\frac{1}{N^{2p-1}}\sum_{n=1}^N(n+1)^{p-2}\left
(\sum_{k=0}^nk\right )^p\asymp
\frac{1}{N^{2p-1}}\sum_{n=1}^N(n+1)^{3p-2}\asymp N^p,$$ it follows
that $\vert \eta _N\vert =\og \left (\frac{1}{N}\right )$.
\end{pf}
\par\smallskip
We remark that the sequence $\{ \eta _n\} $ constructed to prove
part (ii) of Theorem\,\@\ref{og1n} provides us with an example which
shows that the converse of Proposition\,\@\ref{vertetan} is not true
for $1\le p<2$.
\par\bigskip
\section{Rhaly operators acting on derivative-Hardy spaces}\label{RhalySp}
The proof of Theorem\,\@\ref{RSp} will be based on the
Littlewood-Paley characterization of Hardy spaces. Namely, we shall
use the following result which can be found in \cite[Vol. II,
Chapter\,\@XIV,  Theorem\,\@4.\,\@14]{Zy}.
\begin{other}\label{LP-Zy} Suppose that $1<q<\infty $ and let $\{
\lambda _j\} _{j=0}^\infty $ be a sequence of complex numbers. Let
$M$ be a positive number and suppose that
$$\vert \lambda_j\vert \le M,\,\,\,\,\sum_{j=2^k}^{2^{k+1}-1}\vert
\lambda _{j+1}-\lambda _j\vert\,\le M,\quad j, k=0, 1, 2, \dots .$$
There exists a positive constant $A_q$ such that if $f\in H^q$,
$f(z)=\sum_{n=0}^\infty a_nz^n$ ($z\in\mathbb D$), and
$$h(z)=\sum_{n=0}^\infty \lambda_na_nz^n,\quad z\in \mathbb D,$$
then $h\in H^q$ and
$$\Vert h\Vert _{H^q}\,\le \,MA_q\Vert f\Vert _{H^q}.$$
\end{other}\par\medskip
\begin{Pf}{\it Theorem\,\@\ref{RSp}.}
\par The implication (i)\,\,$\Rightarrow$\,\, (iii) is clear.
Indeed,
$$\mathcal R_{(\eta )}\in \mathcal B(S^p, S^q)\,\,\Rightarrow\,\,
\mathcal R_{(\eta )}(1)\in S^q\,\, \Leftrightarrow\,\,F_{(\eta )}\in
S^q.$$
\par\smallskip Let us prove next that (iii)\,\,$\Rightarrow$\,\,
(i). Assume that $F_{(\eta )}\in S^q$, that is, $F_{(\eta )}^\prime
\in H^q$. Take $f\in S^p$, $f(z)=\sum_{n=0}^\infty a_nz^n$ ($z\in
\mathbb D$). Then $f^\prime \in H^1$ and then, using Hardy's
inequality, we see that
$$\sum_{n=0}^\infty \vert a_{n+1}\vert \le \pi \Vert f^\prime \Vert
_{H^1}\le \pi \Vert f^\prime \Vert _{H^p}.$$ This implies that
\begin{equation}\label{sumsp}\sum_{n=0}^\infty \vert a_n\vert \le
\pi \Vert f\Vert _{S^p}.\end{equation} We have
$$\mathcal R_{(\eta )}(f)(z)=\sum_{j=0}^\infty \eta _j\left
(\sum_{k=0}^ja_k\right )z^j,\quad z\in \mathbb D.$$ Hence,
$$\mathcal R_{(\eta )}(f)^\prime (z)=\sum_{j=0}^\infty (j+1)\eta _{j+1}\left
(\sum_{k=0}^{j+1}a_k\right )z^{j},\quad z\in \mathbb D.$$ Set
$$\lambda_j\,=\,\sum_{k=0}^{j+1}a_k,\quad j=0, 1, 2, \dots .$$
Using (\ref{sumsp}), se see that
$$\vert \lambda _j\vert \le \sum_{k=0}^{j+1}\vert a_k\vert \le \sum_{k=0}^{\infty }\vert
a_k\vert \le \pi \Vert f\Vert_{S^p},$$
$$\sum_{j=0}^\infty \vert \lambda_{j+1}-\lambda _j\vert =\sum_{j=0}^\infty \vert a_{j+2}\vert \le \pi \Vert f\Vert_{S^p}.$$
Then Theorem\,\@\ref{LP-Zy} yields that
$$\Vert \mathcal R_{(\eta )}(f)^\prime \Vert _{H^q}\le A_q\pi\Vert
f\Vert _{S^p}\Vert F_{(\eta )}^\prime\Vert _{H^q}.$$ Thus (i)
follows.
\par The implication (ii)\,\,$\Rightarrow$\,\,(i) is obvious. Let us
prove that (iii)\,\,$\Rightarrow$\,\,(ii). So, assume (iii), that
is, $F_{(\eta )}\in S^q$. Then we already know that $\mathcal
R_{(\eta )}\in \mathcal B(S^p, S^q)$.
\par Since $1<q<\infty $, we have that
\begin{equation}\label{partialcon}S_N(F_{(\eta )}^\prime )\to F_{(\eta )}^\prime ,\quad \text{as
$N\to \infty $,\,\,in $H^q$}.\end{equation} For $N=2, 3, 4, \dots $,
let $R_N$ be defined as above, that is, if $f\in S^p$,
$f(z)=\sum_{n=0}^\infty a_nz^n$ ($z\in \mathbb D$), then
$$R_N(f)(z)\,=\,\sum_{j=0}^N\eta _j\left (\sum_{k=0}^ja_k\right
)z^j,\quad z\in \mathbb D.$$ Then
$$\mathcal R_{(\eta )}(f)^\prime (z)-R_N(f)^\prime
(z)\,=\,\sum_{j=N}^\infty (j+1)\eta _{j+1}\left
(\sum_{k=0}^{j+1}a_k\right )z^j.$$ Set

$$
\lambda_j\,=\,\left\{%
\begin{array}{ll}
   0,\,\,\, & \hbox{if $j=0, 1, 2, \dots N$}.\\
   \\
   \sum_{k=0}^{j+1}a_k,\,\,\, & \hbox{if $j=N, N+1, \dots $}.
\end{array}%
\right.
$$
Just as above we have that
$$\vert \lambda_j\vert \le \pi \Vert f\Vert
_{S^p},\,\,\text{for all $j$,\,\, and}\,\,\,\,\sum_{j=0}^\infty
\vert \lambda_{j+1}-\lambda _j\vert \le \pi \Vert f\Vert _{S^p}.$$
Then Theorem\,\@\ref{LP-Zy} implies that
$$\Vert \mathcal R_{(\eta )}(f)^\prime (z)-R_N(f)^\prime
(z)\Vert _{H^q}\le A_q\pi \Vert F_{(\eta )}^\prime -S_{N-1}(F_{(\eta
)}^\prime )\Vert _{H^q}\Vert f\Vert _{S^p}.$$ This and
(\ref{partialcon}) give that $R_N\to \mathcal R_{(\eta )}$ in
$\mathcal B(S^p, S^q)$. Since the operators $R_N$ are finite rank
operators from $S^p$ into $S^q$, it follows that $\mathcal R_{(\eta
)}\in \mathcal K(S^p, S^q)$.
\end{Pf}

\par\medskip
{\bf Data Availability.} All data generated or analyzed during this
study are included in this article and in its bibliography

\par\medskip
{\bf Conflict of interest.} The authors declare that there is no
conflict of interest.

\par
\end{document}